\documentclass[12pt,a4paper]{amsart}
\usepackage[top=25mm, bottom=25mm, left=30mm, right=30mm]{geometry}
\usepackage{amsmath,amssymb,amsthm,eucal,mathrsfs}
\usepackage{verbatim}
\usepackage{bbding}
\usepackage{enumerate}
\usepackage{color}
\usepackage[all]{xy}
\usepackage{mathptmx}
\usepackage[colorlinks=true,citecolor=blue]{hyperref}
\usepackage{xstring}
\renewcommand\MR[1]{%
\StrBefore{#1 }{ }[\USTCtemp]%
\href{http://www.ams.org/mathscinet-getitem?mr=\USTCtemp}{MR\USTCtemp}}
\newcommand{\arXiv}[1]{\href{http://arxiv.org/abs/#1}{arXiv:#1}}

\newtheorem{thm}{Theorem}[section]
\newtheorem{lem}[thm]{Lemma}
\newtheorem{prop}[thm]{Proposition}
\newtheorem{cor}[thm]{Corollary}

\newtheorem{prob}{Problem}
\theoremstyle{definition}
\newtheorem{rem}[thm]{Remark}
\newtheorem{defn}[thm]{Definition}

\newcommand{\eps}{\varepsilon}
\newcommand{\N}{\mathbb{N}}
\newcommand{\F}{\mathcal{F}}

\newcommand{\Z}{{\mathbb{Z}_+}}

\newcommand{\scr}{\mathscr}
\newcommand{\ep}{\varepsilon}

\newcommand{\M}{\overline{\mathcal{M}}}

\DeclareMathOperator{\diam}{diam}

\begin{document}
\title{Recent development of chaos theory in topological dynamics}
\author[J.~Li]{Jian Li}
\date{\today}
\address[J. Li]{Department of Mathematics, Shantou University, Shantou, Guangdong, 515063, P.R. China}
\email{lijian09@mail.ustc.edu.cn}

\author[X.~Ye]{Xiangdong Ye}
\address[X.~Ye]{Wu Wen-Tsun Key Laboratory of Mathematics, USTC, Chinese Academy of Sciences and
School of Mathematics, University of Science and Technology of China,
Hefei, Anhui, 230026, P.R. China}
\email{yexd@ustc.edu.cn}

\begin{abstract}
We give a summary on the recent development of chaos theory in topological dynamics, focusing on
Li-Yorke chaos, Devaney chaos, distributional chaos,
positive topological entropy, weakly mixing sets and so on, and their relationships.
\end{abstract}
\maketitle

\section{Introduction}
Topological dynamics is a branch of the theory of dynamical systems in which qualitative,
asymptotic properties of dynamical systems are studied, %from the viewpoint of general topology,
where a dynamical system is a phase (or state) space $X$ endowed with an evolution map $T$ from $X$ to itself.
In this survey, we require that the phase space $X$ is a compact metric space
and the evolution map $T\colon X\to X$ is continuous.

The mathematical term \emph{chaos} was first introduced by Li and Yorke in 1975~\cite{Li1975},
where the authors established a simple criterion on the existence of chaos for interval maps,
known as ``period three implies chaos''. Since then, the study of chaos theory has played a big role
in dynamical systems, even in nonlinear science.

In common usage, ``chaos'' means ``a state of disorder''.
However, in chaos theory, the term is defined more precisely.
Various alternative, but closely related definitions of chaos have been proposed after Li-Yorke chaos.

Although there is still no definitive, universally accepted mathematical definition of chaos
(in our opinion it is also impossible), most definitions of chaos are based on one of the following aspects:
\begin{enumerate}
  \item complex trajectory behavior of points, such as Li-Yorke chaos and distributional chaos;
  \item sensitivity dependence on initial conditions, such as Devaney  chaos,
   Auslander-Yorke chaos and Li-Yorke sensitivity;
  \item fast growth of different orbits of length $n$, such as having positive topological entropy, or its variants;
  \item strong recurrence property, such as weakly mixing property and  weakly mixing sets.
\end{enumerate}

Since there are so many papers dealing with the chaos theory in topological dynamics,
we are not able to give a survey on all of them. We only can select those which we are familiar with,
and are closely related to our interest, knowledge and ability.
As we said before it is impossible to give a universe definition of chaos which covers all the
features of the complex behaviors of a dynamical system, it is thus important to know the relationships among
the various definitions of chaos. So we will focus on Li-Yorke chaos, Devaney chaos, distributional chaos,
positive topological entropy, weakly mixing sets, sensitivity and so on, and their relationships in this survey.
See~\cite{Blanchard2009,Oprocha2014} for related surveys. %discussion on various kinds of chaos.

\medskip

Now we outline the development of chaos theory for general topological dynamical systems briefly.
The notions of topological entropy and weak mixing were introduced by Adler et al \cite{Adler1965}
in 1965 and Furstenberg \cite{Furstenberg1967} in 1967 respectively. After the mathematical term of
chaos by Li-Yorke appearing in 1975, Devaney \cite{Devaney1989}
defined a kind of chaos, known as Devaney chaos today,
in 1989 based on the notion of sensitivity introduced by Guckenheimer \cite{Guckenheimer1979}.
The implication among them has attracted a lot of
attention. In 1991 Iwanik \cite{Iwanik1991} showed that weak mixing implies Li-Yorke chaos.
A breakthrough concerning the relationships among positive entropy, Li-Yorke chaos and Devaney chaos came in 2002.
In that year, it was shown by Huang and Ye \cite{Huang2002} that Devaney
chaos implies Li-Yorke one by proving that a non-periodic transitive
system with a periodic point is Li-Yorke chaotic;
Blanchard, Glasner, Kolyada and Maass \cite{Blanchard2002} proved
that positive entropy also implies Li-Yorke chaos (we remark that the
authors obtained this result using ergodic method, and for a
combinatorial proof see \cite{Kerr2007}). Moreover, the result also holds
for sofic group actions by Kerr and Li \cite[Corollary 8.4]{Kerr2013}.

In 1991, Xiong and Yang~\cite{Xiong1991} showed that in a weakly mixing system
there are considerably many points in the domain whose orbits display
highly erratic time dependence.
It is known that a dynamical system with positive entropy may not contain any weakly mixing subsystem.
Capturing a common feature of positive entropy and weak mixing, Blanchard and Huang \cite{Blanchard2008} in 2008
defined the notion of weakly mixing set and showed that positive entropy implies the existence of weakly mixing sets
which also implies Li-Yorke chaos.
A further discussion along the line will be appeared in a forthcoming paper
by Huang, Li, Ye and Zhou \cite{Huang2014d}.

Distributional chaos was introduced in 1994 by
Schweizer and Sm{\'{\i}}tal \cite{Schweizer1994}, and there are at
least three versions of distributional chaos in the literature (DC1,
DC2 and DC3). It is known that positive entropy does not imply DC1 chaos \cite{Pikula2007} and
Sm\'{\i}tal conjectured that positive entropy implies DC2 chaos.
Observing that DC2 chaos is equivalent to so called mean Li-Yorke chaos,
recently Downarowicz \cite{Downarowicz2011} proved that positive entropy indeed implies DC2 chaos.
An alternative proof can be found in \cite{Huang2014} by Huang, Li and Ye. We remark that both proofs
use ergodic theory heavily and there is no combinatorial proof at this moment.

This survey will be organized as follows. In Section \ref{sec:pre} we provide basic definitions in topological dynamics.
Li-Yorke chaos, sensitivity and chaos in transitive systems will be discussed in Sections \ref{sec:LY-Chaos}-\ref{section4}.
In Section \ref{section5} we review the results on distributional chaos.
In the following two sections we focus on  weakly mixing sets
and chaos in the induced spaces.

\section{Preliminaries}\label{sec:pre}
In this section, we provide some basic notations, definitions and results which will be
used later in this survey.

Denote by $\N$ ($\Z$, $\mathbb{Z}$, respectively) the set of all positive integers
(non-negative integers, integers, respectively).
The cardinality of a set $A$ is usually denoted by $|A|$.

Let $X$ be a compact metric space.
A subset $A$ of $X$ is called a \emph{perfect set} if it is a closed set with no isolated points;
a \emph{Cantor set} if it is homeomorphic to  the standard  middle third  Cantor set;
a \emph{Mycielski set} if it is a union of countably many Cantor sets.
For convenience we restate here a
version of Mycielski's theorem (\cite[Theorem~1]{Mycielski1964}) which we shall use.

\begin{thm}[Mycielski Theorem]\label{thm:Mycielski-thm}
Let $X$ be a perfect compact metric space.
If $R$ is a dense $G_\delta$ subset of $X^n$, then there exists a
dense Mycielski set $K\subset X$ such that for any distinct $n$ points $x_1,x_2,\dotsc,x_n\in K$,
the tuple $(x_1,x_2,\dotsc,x_n)$ is in $R$.
\end{thm}

\subsection{Topological dynamics}
By a \emph{(topological) dynamical system}, we mean a pair $(X,T)$, where $X$ is a compact metric space
and $T\colon X\to X$ is a continuous map. The metric on $X$ is usually denoted by $d$.
One can think of $X$ as a state space for some system, and $T$ as
the evolution law of some discrete  autonomous  dynamics on $X$:
if $x$ is a point in $X$, denoting the current state of a system,
then $Tx$ can be interpreted as the state of the same system after one
unit of time has elapsed.
For every non-negative integer $n$, we can define the iterates $T^n\colon X\to X$ as
$T^0=id $ the identity map on $X$  and $T^{n+1}=T^n\circ T$.
One of the main topics of study in dynamical systems is the asymptotic
behaviour of $T^n$ as $n \to\infty$.

Note that we always assume that the state space $X$ is not empty.
If the state space $X$ contains only one point, then we say that the dynamical system on $X$ is \emph{trivial},
because in this case the unique map on $X$ is the identity map.

For any $n\geq 2$, the \emph{$n$-th fold product} of $(X,T)$ is denoted by $(X^{n}, T^{(n)})$,
where $X^{n}=X\times X\times \dotsb \times X$ ($n$-times) and
$T^{(n)}=T\times T\times \dotsb \times T$ ($n$-times).
We set the diagonal of $X^{n}$ as  $\Delta_n=\{(x, x,\dotsc, x)\in X^{n}\colon x\in X\}$,
and set $\Delta^{(n)}=\{(x_1, x_2,\dotsc,x_n)\in X^{n}\colon
\text{there exist }1\leq i<j\leq n\text{ such that }x_i=x_j\}$.

The \emph{orbit} of a point $x$ in $X$ is the set $Orb(x,T)=\{x,Tx,T^2x,\dotsc\}$;
the $\omega$-limit set of $x$, denoted by $\omega(x,T)$, is the limit set of the orbit of $x$,
that is
\[\omega(x,T)=\bigcap_{n=1}^\infty\overline{\{T^ix\colon i\geq n\}}.\]
A point $x\in X$ is called \emph{recurrent} if there exists an increasing sequence $\{p_i\}$ in $\N$
such that $\lim_{i\to\infty} T^{p_i}x=x$. Clearly, $x$ is recurrent if and only if
$x\in\omega(x,T)$.

If $Y$ is a non-empty closed invariant (i.e.\ $TY\subset Y$) subset of $X$,
then $(Y,T)$ is also a  dynamical system,
we call it as a \emph{subsystem} of $(X,T)$.
A dynamical system $(X,T)$ is called \emph{minimal} if it contains no proper subsystems.
A subset $A$ of $X$ is \emph{minimal} if $(A,T)$ forms a minimal subsystem of $(X,T)$.
It is easy to see that a non-empty closed invariant set $A\subset X$ is minimal if and only if
the orbit of every point of $A$ is dense in $A$.
A point $x\in X$ is called \emph{minimal} or \emph{almost periodic}
if it belongs to a minimal set.

A dynamical system $(X,T)$ is called \emph{transitive} if for every two non-empty open subsets $U$ and $V$
of $X$ there is a positive integer $n$ such that $T^n U\cap V\neq\emptyset$;
\emph{totally transitive} if $(X,T^n)$ is transitive for all $n\in\N$;
\emph{weakly mixing} if the product system $(X\times X,T\times T)$ is transitive;
\emph{strongly mixing} if for every two non-empty open subsets $U$ and $V$
of $X$ there is $N>0$ such that $T^n U\cap V\neq\emptyset$ for all $n\geq N$.
Any point with dense orbit is called a \emph{transitive point}.
In a transitive system the set of all transitive points is a dense $G_\delta$ subset of $X$
and we denote it by $Trans(X,T)$. For more details related to transitivity see~\cite{Kolyada1997}.

Let $(X,T)$ be a  dynamical system.
A pair $(x,y)$ of points in $X$ is called
\emph{asymptotic} if
  $\lim\limits_{n\to\infty} d(T^nx,T^ny)=0$;
\emph{proximal} if
  $\liminf\limits_{n\to\infty} d(T^nx,T^ny)=0$;
\emph{distal} if
 $\liminf\limits_{n\to\infty} d(T^nx,T^ny)>0$.
The system $(X,T)$ is called \emph{proximal} if any two points in $X$ form a  proximal pair;
\emph{distal} if any two distinct points in $X$ form a distal pair.

Let $(X,T)$ and $(Y,S)$ be two dynamical systems.
If there is a continuous surjection $\pi\colon X\to Y$ which intertwines the actions (i.e., $\pi\circ T=S\circ \pi$),
then we say that $\pi$ is a \emph{factor map},
 $(Y,S)$ is a \emph{factor} of $(X,T)$ or $(X,T)$ is an \emph{extension} of $(Y,S)$.
The factor map $\pi$ is  \emph{almost one-to-one} if there exists a residual subset $G$ of $X$ such that
$\pi^{-1}(\pi(x))=\{x\}$ for any $x\in G$.

In 1965, Adler, Konheim and McAndrew introduced topological entropy in topological dynamics~\cite{Adler1965}.
Let $C^o_X$ be the set of finite open covers of $X$.
Given two open covers $\mathcal U$ and $\mathcal V$, their \emph{join} is the cover
\[\{U\cap V\colon U\in\mathcal{U},V\in\mathcal{V}\}.\]
We define $N(\mathscr U)$ as the minimum cardinality of subcovers of $\mathcal U$.
The \emph{topological entropy} of $T$ with respect to $\mathcal{U}$ is
\begin{equation*}%\label{TE}
h_{\text{top}}(T,\mathcal{U})=\lim_{N\to\infty}{\frac1N}\log N\bigg(\bigvee_{i=0}^{N-1}T^{-i}\mathcal{U}\bigg).
\end{equation*}
The \emph{topological entropy} of $(X,T)$ is defined by
\[h_{\text{top}}(T)=\sup_{\mathcal{U}\in C^o_X}h_{\text{top}}(T,\mathcal{U}).\]
We refer the reader to \cite{Walters1982,Glasner2003} for more information on topological entropy,
and~\cite{Glasner2009} for a survey of local entropy theory.

\subsection{Furstenberg family}
The idea of using families of subsets of $\Z$ when dealing with dynamical properties of maps was
first used by Gottschalk and Hedlund~\cite{Gottschalk1955},
and then further developed by Furstenberg~\cite{Furstenberg1981}.
For a systematic study and recent results,
see~\cite{Akin1997,Huang2004a,Huang2004,Li2011b,Shao2006,Xiong2007}.
Here we recall some basic facts related to Furstenberg  families.

A \emph{Furstenberg family} (or just \emph{family}) $\F$ is a collection of subsets of $\Z$ which is upwards hereditary,
that is, if $F_1\in\F$ and $F_1\subset F_2$ then $F_2\in\F$.
A family $\F$ is \emph{proper} if $\N\in\F$ and $\emptyset\not\in\F$.
If $\F_1$ and $\F_2$ are families, then we define
\[\F_1\cdot \F_2  = \{F_1 \cap F_2\colon F_1\in\F_1,F_2\in\F_2\}.\]
If $\F$ is a family, the \emph{dual family} of $\F$, denoted by $\kappa\F$, is the family
\[\{F\subset\N\colon F\cap F'\neq\emptyset, \forall F'\in\F\}.\]
We denote by $\F_{inf}$ the family of infinite subsets of $\N$.
The dual family of $\F_{inf}$ is just the collection of co-finite subsets of $\Z$,
denoted by $\F_{cf}$.
A family $\F$ is called \emph{full} if it is proper and $\F\cdot\kappa\F\subset\F_{inf}$.

For a sequence $\{p_i\}_{i=1}^{\infty}$ in $\N$, define the finite sums of $\{p_i\}_{i=1}^{\infty}$ as
\[FS\{p_i\}_{i=1}^{\infty}=\biggl\{ \sum_{i\in \alpha}p_i\colon
\alpha \textrm{ is a non-empty finite subset of }\N\,\biggr\}.\]
A subset $F$ of $\Z$ is called an {\em IP-set}
if there exists a sequence $\{p_i\}_{i=1}^{\infty}$ in $\N$
such that $ FS\{p_i\}_{i=1}^{\infty}\subset F$.
We denote by $\F_{ip}$ the family of all IP-sets.

A subset $F$ of $\Z$ is called \emph{thick}
if it contains arbitrarily long runs of positive integers, i.e., for every $n\in \N$ there exists
some $a_n \in\Z$ such that $\{a_n,a_n+1,\dotsc,a_n+n\}\subset F$;
\emph{syndetic} if it has bounded gaps, i.e., there is $N\in\N$ such that $[n,n + N] \cap F\neq\emptyset$
for every  $n \in\Z$. The families of all thick sets and syndetic sets are denoted by
$\F_t$ and $\F_s$, respectively. It is easy to see that $\kappa\F_s = \F_t$.

Let $(X, T)$ be a dynamical system.
For $x \in X$ and a non-empty subset $U$ of $X$, we define the \emph{entering time set} of $x$ into $U$ as
\[ N(x,U) =\{n \in\mathbb{Z}_+\colon T^nx \in  U  \}.\]
If $U$ is a neighborhood of $x$, then we usually call $N(x,U)$ the \emph{return time set} of $x$ into $U$.
Clearly, a point $x$ is recurrent if and only if
for every open neighborhood $U$ of $x$, the return time
set $N(x,U)$ is infinite.
In general, we can define recurrence with respect to a Furstenberg family.
\begin{defn}
Let $(X, T)$ be a dynamical system and $\F$ be a Furstenberg family.
A point $x\in X$ is said to be \emph{$\F$-recurrent} if for every open neighborhood $U$ of $x$,
$N(x,U)\in\F$.
\end{defn}

It is well known that the following
lemma  holds (see, e.g., \cite{Akin1997,Furstenberg1981}).
\begin{lem}
Let $(X,T)$ be a dynamical system and $x\in X$. Then
\begin{enumerate}
  \item $x$ is a minimal point if and only if it is an $\F_s$-recurrent point;
  \item $x$ is a recurrent point if and only if it is an $\F_{ip}$-recurrent point.
\end{enumerate}
\end{lem}

For two non-empty subsets $U$ and $V$ of $X$, we define the \emph{hitting time set} of $U$ and $V$ as
\[ N(U, V ) =\{n \in\mathbb{Z}_+\colon T^nU \cap V \neq\emptyset\}.\]
If $U$ is a non-empty open subset of $X$, then we usually call $N(U,U)$ the \emph{return time set} of $U$.
Clearly, a dynamical system $(X,T)$ is transitive (resp.\ strongly mixing)
if and only if for any two non-empty open subsets $U$, $V$ of $X$,
the hitting time set $N(U,V)$ is infinite (resp.\ co-finite).

\begin{defn}
Let $(X, T)$ be a dynamical system and $\F$ be a family. We say that $(X, T)$  is
\begin{enumerate}
  \item \emph{$\F$-transitive} if for any two non-empty open subsets $U$, $V$ of $X$, $N(U,V)\in\F$;
  \item \emph{$\F$-mixing} if the product system $(X\times X,T\times T)$ is $\F$-transitive.
\end{enumerate}
\end{defn}

\begin{thm}[\cite{Furstenberg1967}]
Let $(X,T)$ be a dynamical system.
Then $(X,T)$ is weakly mixing if and only if
it is $\F_t$-transitive.
\end{thm}

\section{Li-Yorke chaos} \label{sec:LY-Chaos}
The mathematical terminology ``chaos'' was first introduced in 1975 by Li and Yorke to describe the complex
behavior of trajectories. It turns out that it is the common feature of all known definitions of chaos.
\subsection{Li-Yorke chaos and its relation with topological entropy}
Following the idea in~\cite{Li1975},
we usually define the Li-Yorke chaos as follows.
\begin{defn}
Let $(X,T)$ be a dynamical system.
A pair $(x,y)\in X\times X$ is called \emph{scrambled } if
\[ \liminf_{n\to\infty} d(T^n x,T^ny)=0 \quad \text{and}\quad \limsup_{n\to\infty} d(T^n x,T^ny)>0,\]
that is $(x,y)$ is proximal but not asymptotic.
A subset $C$ of $X$ is called \emph{scrambled} if any two distinct points $x,y\in C$ form a scrambled pair.
The dynamical system $(X,T)$ is called \emph{Li-Yorke chaotic} if there is an uncountable  scrambled set in $X$.
\end{defn}

It is worth to notice that the terminology ``scrambled set'' was introduced in 1983 by Sm\'{\i}tal in~\cite{Smital1983}.
Note that we should assume that a scrambled set contains at least two points.
We will keep this convention throughout this paper.

In~\cite{Li1975}, Li and Yorke showed that
\begin{thm}\label{thm:LY3}%[\cite{Li1975}]
If a continuous map $f\colon[0,1]\to [0,1]$ has a periodic point of period $3$,
then it is Li-Yorke chaotic.
\end{thm}

\begin{defn}
Let $(X,T)$ be a dynamical system.
For a given positive number $\delta>0$, a pair $(x,y)\in X\times X$
is called \emph{$\delta$-scrambled}, if
\[ \liminf_{n\to\infty} d(T^n x,T^ny)=0 \quad
\text{and} \quad \limsup_{n\to\infty} d(T^n x,T^ny)>\delta.\]
A subset $C$ of $X$ is \emph{$\delta$-scrambled} if any two distinct points $x,y$ in $C$
 form a $\delta$-scrambled pair.
The dynamical system $(X,T)$ is called \emph{Li-Yorke $\delta$-chaotic},
if there exists an uncountable $\delta$-scrambled set in $X$.

Note that the Auslander-Floyd system~\cite{Auslander1988} is Li-Yorke chaotic,
but for any $\delta>0$ there is no uncountable  $\delta$-scrambled sets.
\end{defn}

In~\cite{Jankova1986} Jankov{\'a} and Sm{\'{\i}}tal showed that
if a continuous map $f:[0,1]\to[0,1]$ has positive topological entropy then
there exists a perfect $\delta$-scrambled set for some $\delta>0$.
A natural question is that whether there exists a Li-Yorke chaotic map with
zero topological entropy.
This was shown, independently, by Xiong~\cite{Xiong1986} and Sm{\'{\i}}tal~\cite{Smital1986} in 1986.

\begin{thm} \label{thm:LY-chaos-0entropy}
There exists a continuous map $f\colon[0,1]\to [0,1]$,
which is Li-Yorke chaotic but has zero topological entropy.
\end{thm}

In 1991, using Mycielski Theorem~\ref{thm:Mycielski-thm}, Iwanik~\cite{Iwanik1991} showed
that every weakly mixing system is Li-Yorke chaotic.
\begin{thm}[\cite{Iwanik1991}]
If a non-trivial dynamical system $(X,T)$ is weakly mixing, then there exists a dense
Mycielski $\delta$-scrambled subset of $X$ for some $\delta>0$.
\end{thm}

As already mentioned  above, for interval maps positive topological entropy implies Li-Yorke chaos.
It used to be a long-standing open problem whether this also holds for general topological dynamical systems.
In 2002, Blanchard, Glasner, Kolyada and Maass gave a positive answer by using ergodic theory method.

\begin{thm}[\cite{Blanchard2002}]\label{thm:BGKM}
If a dynamical system $(X,T)$ has positive topological entropy, then
there exists a Mycielski  $\delta$-scrambled set for some $\delta>0$.
\end{thm}

In fact we know from the proof of Theorem~\ref{thm:BGKM}
that if $(X,T)$ has an ergodic invariant measure which is not measurable distal then
same conclusion holds.
In 2007, Kerr and Li~\cite{Kerr2007} gave a new proof of Theorem~\ref{thm:BGKM} by using combinatorial method.
First we recall the definition of IE-tuples.

\begin{defn}\label{def:ind-sets}
Let $(X,T)$ be a topological dynamical system and $k\geq 2$.
For a tuple $\tilde A=(A_1,\dotsc,A_k)$ of subsets of $X$,
we say that a subset $J$ of $\Z$ is an \emph{independence set} for $\tilde A$
if for any non-empty finite subset $I$ of $J$, we have
\[\bigcap_{i\in I}T^{-i} A_{s(i)}\neq\emptyset.\]
for any $s\in\{1,\dotsc,k\}^I$.

A tuple $\tilde x=(x_1,\dotsc,x_k)\in X^k$ is called an \emph{IE-tuple}
if for every product neighborhood $U_1\times \dotsb\times U_k$ of $\tilde x$
the tuple $(U_1,\dotsc,U_k)$ has an independence set of positive density.
\end{defn}

The following theorem characterizes positive topological entropy and has many applications.
It was first proved by Huang and Ye using the notion of interpolating sets \cite{Huang2006},
and we state it here using the notion of independence by Kerr and Li \cite{Kerr2007}.
Recall for a dynamical system $(X,T)$, a tuple $(x_1,\ldots,x_k)\in X^k$
is an \emph{entropy tuple} if $x_i\neq x_j$ for $i\neq j$ and
for any disjoint closed neighborhoods $V_i$ of $x_i$,
the open cover $\{V_1^c, \ldots, V_k^c\}$ has positive entropy.
When $k=2$ we call it an \emph{entropy pair}.
A subset $A$ of $X$ is an \emph{entropy set} if
any tuple of points in $A$ with pairwise different coordinates is an entropy tuple.

\begin{thm}[\cite{Huang2006,Kerr2007}]
Let $(X,T)$ be a dynamical system.
Then a tuple on $X$ is an entropy tuple if and only if it is a non-diagonal IE-tuple.
In particular, the system $(X,T)$ has zero topological entropy if and only if
every IE-pair is diagonal (i.e. all of its entries are equal).
\end{thm}

By developing some deep combinatorial tools Kerr and Li showed

\begin{thm}[\cite{Kerr2007}]\label{thm:IE-chaos}
Let $(X,T)$ be a dynamical system.
Suppose that $k\geq 2$ and $\tilde x=(x_1,\dotsc,x_k)$ is a non-diagonal IE-tuple.
For each  $1\leq j\leq k$, let $A_j$ be a neighborhood of $x_j$.
Then there exists a Cantor set $Z_j\subset A_j$ for each $j=1,\dotsc,k$ such that the following hold:
\begin{enumerate}
  \item every tuple of finite points in $Z:=\bigcup_j Z_j$ is an IE-tuple;
  \item for all $m\in\N$, distinct $y_1,\dotsc,y_m\in Z$ and $y_1',\dotsc,y_m'\in Z$ one has
  \[\liminf_{n\to\infty} \max_{1\leq i\leq m} d(T^n y_i,y_i')=0.\]
\end{enumerate}
In particular, $Z$ is $\delta$-scrambled for some $\delta>0$.
\end{thm}

Recently, Kerr and Li showed that Theorem~\ref{thm:IE-chaos} also holds for sofic group actions
(see \cite[Corollary 8.4]{Kerr2013}).

\subsection{Completely scrambled systems and invariant scrambled sets}
\begin{defn}
We say that a dynamical system $(X,T)$ is
\emph{completely scrambled} if the whole space $X$ is scrambled.
\end{defn}
It should be noticed that for any $\delta>0$,
the whole space can not be $\delta$-scrambled~\cite{Blanchard2008a}.
Recall that a dynamical system $(X,T)$ is proximal if any two points in $X$ form a proximal pair.
Clearly, every completely scrambled system is proximal.
We have the following characterization of proximal systems.
\begin{thm}[\cite{Huang2011,Akin2003}]
A dynamical system $(X,T)$ is proximal if and only if it has a fixed point which is the
unique minimal point in $X$.
\end{thm}

In 1997, Mai~\cite{Mai1997} showed that there are some completely scrambled systems on non-compact spaces.
\begin{thm}
Let $X$ be a metric space uniformly homeomorphic to the $n$-dimensional open cube $I^n=(0,1)^n$, $n\geq 2$.
Then there exists a homeomorphism $f\colon X\to X$ such that the whole space $X$
is a scrambled set of $X$.
\end{thm}

In 2001, Huang and Ye~\cite{Huang2001} showed that on some compact spaces there are also some
completely scrambled systems.
\begin{thm}\label{thm:completely-scrambled-systems}
There are ``many'' compacta admitting completely scrambled homeomorphisms,
which include some countable compacta, the Cantor set and continua of arbitrary dimension.
\end{thm}

Huang and Ye also mentioned in~\cite{Huang2001} that an example of completely scrambled
transitive homeomorphism is a consequence of construction of uniformly rigid proximal systems
by Katznelson and Weiss~\cite{Katznelson1981}.
Recall that a dynamical system $(X,T)$ is \emph{uniformly rigid} if
\[\liminf_{n\to\infty}\sup_{x\in X} d(T^n(x),x)=0.\]
Later in~\cite{Huang2002},
they showed that every almost equicontinuous but not minimal system has a completely scrambled factor.
These examples are not weakly mixing,
so the existence of completely scrambled weakly mixing homeomorphism is left open in~\cite{Huang2001}.
Recently, For\'ys \textit{et al.} in~\cite{Forys2013} constructed two kinds of completely scrambled systems
which are weakly mixing, proximal and uniformly rigid. The first possible approach is
derived from results of Akin and Glasner \cite{Akin2001} by a combination of abstract arguments.
The second method is obtained by modifying the construction of Katznelson and Weiss
from~\cite{Katznelson1981}.
More precisely, we have the following result.
\begin{thm}[\cite{Forys2013}] \label{thm:wm-uniformly-rigid}
There are completely scrambled systems which are weakly mixing, proximal and uniformly rigid.
\end{thm}

On the other hand,  Blanchard, Host and Ruette~\cite{Blanchard2002a} proved that any  positive
topological entropy system can not be  completely scrambled by showing that
there are ``many'' non-diagonal asymptotic pairs in
any dynamical system with positive topological entropy.
\begin{thm}[\cite{Blanchard2002a}] \label{thm:asy-pairs}
Let $(X,T)$ be a dynamical system and $\mu$ be an ergodic invariant
measure on $(X,T)$ with positive entropy. Then for $\mu$-a.e.\ $x\in X$
there exists a point $y\in X\setminus\{x\}$ such that
$(x,y)$ is asymptotic.
\end{thm}

Recently in \cite{Huang2014c}, Huang, Xu and Yi generalized Theorem~\ref{thm:asy-pairs}
to positive entropy $G$-systems for certain countable, discrete, infinite
left-orderable amenable groups $G$.
We remark that the following question remains open.
\begin{prob}
Let $T$ be a homeomorphism on a compact metric space $X$.
If for any two distinct points $x,y\in X$, $(x,y)$ is either Li-Yorke scrambled for $T$ or Li-Yorke scrambled for $T^{-1}$,
does $T$ have zero topological entropy?
\end{prob}

Let $S$ be a scrambled set of a dynamical system $(X,T)$.
It is easy to see that for any $n\geq 1$, $T^n S$ is also a scrambled set of $(X,T)$.
It is interesting to consider that whether a scrambled set may be invariant under $T$.
Since the space $X$ is compact, if $(x,f(x))$ is scrambled for some $x\in X$,
then there is a fixed point in $X$.
It is shown in~\cite{Du2005} that
\begin{thm}
Let $f:[0,1]\to[0,1]$ be a continuous map.
Then $f$ has positive topological entropy if and only if
$f^n$ has an uncountable invariant scrambled set for some $n>0$.
\end{thm}

In 2009, Yuan and L\"u proved that
\begin{thm}[\cite{Yuan2009}]
Let $(X,T)$ be a non-trivial transitive system.
If $(X,T)$ has a fixed point, then there exists a dense  Mycielski subset $K$ of $X$
such that $K$ is an invariant scrambled set.
\end{thm}

In 2010, Balibrea, Guirao and Oprocha studied invariant $\delta$-scrambled sets and showed that
\begin{thm}[\cite{Balibrea2010}]
Let $(X,T)$ be a non-trivial strongly mixing system.
If $(X,T)$ has a fixed point, then there exist $\delta>0$
and a dense Mycielski subset $S$ of $X$ such that $S$ is an invariant $\delta$-scrambled set.
\end{thm}

They also conjectured in \cite{Balibrea2010} that there exists a weakly mixing system which has a fixed point
but without invariant $\delta$-scrambled sets.
The authors in~\cite{Forys2013} found that the existence of invariant $\delta$-scrambled sets
is relative to the property of uniform  rigidity.
It is shown in~\cite{Forys2013} that a necessary condition for a dynamical system
possessing invariant $\delta$-scrambled sets for some $\delta>0$  is not uniformly rigid,
and this condition (with a fixed point) is also sufficient for transitive systems.

\begin{thm}[\cite{Forys2013}] \label{thm:not-uniformly-rigid}
Let $(X,T)$ be a non-trivial transitive system. Then $(X,T)$ contains
a dense Mycielski invariant $\delta$-scrambled set for some $\delta>0$
if and only if it has a fixed point and is not uniformly rigid.
\end{thm}

Combining Theorems~\ref{thm:wm-uniformly-rigid} and~\ref{thm:not-uniformly-rigid},
the above mentioned conjecture from~\cite{Balibrea2010} has an affirmative answer.
This is because
by Theorem~\ref{thm:wm-uniformly-rigid} there exist weakly mixing, proximal and uniformly rigid systems which have a fixed point, but by Theorem~\ref{thm:not-uniformly-rigid} they do not have any invariant $\delta$-scrambled set.

As far as we know, whenever a dynamical system has been shown to be Li-Yorke chaotic the proof
implies the existence of a Cantor or Mycielski scrambled set.
This naturally arises the following problem:

\begin{prob}\label{prob:cantro-scramble}
If a dynamical system is Li-Yorke chaotic,
dose there exist a Cantor scrambled set?
\end{prob}

Although we do not known the answer of Problem~\ref{prob:cantro-scramble},
there are severe restrictions on the Li-Yorke chaotic dynamical systems without a Cantor scrambled set,
if they exist (see~\cite{Blanchard2008a}). On the other hand, we have
\begin{thm}[\cite{Blanchard2008a}]
If a dynamical system $(X,T)$ is Li-Yorke $\delta$-chaotic for some $\delta>0$,
then it does have some Cantor $\delta$-scrambled set.
\end{thm}

The key point is that the collection of $\delta$-scrambled pairs is a $G_\delta$ subset of $X\times X$,
then one can apply the Mycielski Theorem to get a Cantor $\delta$-scrambled set.
But in general the collection of scrambled pairs is a Bore set of $X\times X$ but may be not $G_\delta$.
By~\cite[Example 3.6]{Snoha1990} there exists a dynamical system such that
the collection of scrambled pairs is dense in $X\times X$ but not residual.

\section{Sensitive dependence on initial conditions}
We say that a dynamical system $(X,T)$ has \emph{sensitive dependence on initial conditions}
(or just \emph{sensitive})
if there exists some $\delta>0$ such that for each $x\in X$ and each $\ep>0$
there is $y\in X$ with $d(x,y)<\ep$ and $n\in\N$
such that $d(T^nx,T^ny)>\delta$.
The initial idea goes back at least to Lorentz \cite{Lorenz1963}
and the phrase---sensitive dependence on initial conditions---was
used by Ruelle \cite{Ruelle1978} to indicate some exponential rate of divergence
of orbits of nearby points.
As far as we know the first to formulate the sensitivity
was Guckenheimer, \cite{Guckenheimer1979}, in his study on maps of the interval
(he required the condition to hold for a set of positive Lebesgue measure).
The precise expression of sensitivity in the above form was introduced by Auslander and Yorke in \cite{Auslander1980}.

\subsection{Equicontinuity and Sensitivity}
The opposite side of sensitivity is the notion of equicontinuity.
\begin{defn}
A dynamical system $(X,T)$ is called \emph{equicontinuous}
if for every $\varepsilon>0$ there exists some $\delta>0$ such that
whenever $x,y\in X$ with $d(x,y)<\delta$, $d(T^nx,T^ny)<\varepsilon$ for $n=0,1,2,\dotsc$,
that is the family of maps $\{T^n\colon n\in\Z\}$ is uniformly equicontinuous.
\end{defn}

Equicontinuous systems have simple dynamical behaviors.
It is well known that a dynamical system $(X,T)$ with $T$ being surjective is equicontinuous if and only if
there exists a compatible metric $\rho$ on $X$ such that $T$ acts on $X$ as an isometry, i.e.,
$\rho(Tx,Ty)=\rho(x,y)$ for any $x,y\in X$.
See~\cite{Mai2003} for the structure of equicontinuous systems.

We have the following dichotomy result for minimal systems.

\begin{thm}[\cite{Auslander1980}]
A minimal system is either equicontinuous or sensitive.
\end{thm}

Equicontinuity can be localized in an obvious way.

\begin{defn}
Let $(X,T)$ be a dynamical system.
A point $x\in X$ is called an \emph{equicontinuous point}
if for any $\ep>0$, there exists some $\delta>0$ such that
$d(x,y)<\delta$ implies $d(T^nx,T^ny)<\ep$ for all $n\in\N$.
A transitive system is called \emph{almost equicontinuous}
if there exists some equicontinuous point.
\end{defn}

We have the following dichotomy result for transitive systems.

\begin{thm}[\cite{Akin1996,Glasner1993}] \label{thm:AY-thm}
Let $(X,T)$ be a transitive system.
Then either
\begin{enumerate}
  \item $(X,T)$ is almost equicontinuous, in this case
  the collection of equicontinuous points coincides with the collection of transitive points; or
  \item $(X,T)$ is sensitive.
\end{enumerate}
\end{thm}

It is interesting that almost equicontinuity is closely related to the uniform rigid property
which was introduced by Glasner and Maon in~\cite{Glasner1989} as a topological
analogue of rigidity in ergodic theory.

\begin{thm}[\cite{Glasner1993}] \label{thm:almost-equi-rigid}
Let $(X,T)$ be a transitive system.
Then it is uniformly rigid if and only if it is a factor of an almost equicontinuous system.
In particular, every almost equicontinuous system is uniformly rigid.
\end{thm}

It is shown in~\cite{Glasner1989}  that every uniformly rigid system has zero topological entropy.
Then by Theorem~\ref{thm:almost-equi-rigid} every almost equicontinuous system also has zero topological entropy.

\subsection{\texorpdfstring{$n$}{n}-sensitivity and sensitive sets}
Among other things, Xiong~\cite{Xiong2005} introduced a new notion called $n$-sensitivity,
which says roughly that in each non-empty open subset there are $n$ distinct points whose
trajectories are apart from (at least for one common moment) a given positive constant pairwise.
\begin{defn}
Let $(X,T)$ be a dynamical system and $n\geq 2$.
The system $(X,T)$ is called \emph{$n$-sensitive},
if there exists some $\delta>0$ such that for any $x\in X$ and $\ep>0$
there are $x_1,x_2,\dotsc,x_n\in B(x,\ep)$ and $k\in\N$ satisfying
\[\min_{1\leq i<j\leq n} d(T^kx_i,T^kx_j) >\delta.\]
\end{defn}

\begin{prop}[\cite{Xiong2005}]
If a dynamical system $(X,T)$ is weakly mixing, then it is $n$-sensitive for all $n\geq 2$.
Moreover, if a dynamical system $(X,T)$ on a locally connected space $X$ is sensitive,
then it is $n$-sensitive for all $n\geq 2$.
\end{prop}

In~\cite{Shao2008}, Shao, Ye and Zhang studied the properties of $n$-sensitivity for minimal systems,
and showed that $n$-sensitivity and $(n+1)$-sensitivity are essentially different.

\begin{thm}[\cite{Shao2008}]
For every $n\geq 2$, there exists a minimal system $(X,T)$
which is $n$-sensitive but not $(n+1)$-sensitive.
\end{thm}

Recently, using ideas and results from local entropy theory, Ye and Zhang~\cite{Ye2008}
and Huang, Lu and Ye~\cite{Huang2011} developed a theory of sensitive sets,
which measures the ``degree'' of sensitivity both in the topological and
the measure-theoretical setting.

\begin{defn}
Let $(X,T)$ be a dynamical system.
A subset $A$ of $X$ is \emph{sensitive} if for any $n\geq 2$, any $n$ distinct points $x_1,x_2,\dotsc,x_n$ in $A$,
any neighborhood $U_i$ of $x_i$, $i=1,2,\dotsc,n$, and any
non-empty open subset $U$ of $X$ there exist $k\in\N$ and $y_i\in U$ such that $T^k(y_i)\in U_i$ for $i=1,2,\dotsc,n$.
\end{defn}

It is shown in~\cite{Ye2008} that a transitive system is $n$-sensitive if and only if there
exists a sensitive set with cardinality $n$. Moreover, a dynamical system is weakly mixing
if and only if the whole space $X$ is sensitive.

\begin{thm}[\cite{Ye2008}]
If a dynamical system is transitive, then every entropy set is also a sensitive set.
This implies that if a transitive system has positive topological entropy, then there
exists an uncountable sensitive set.
\end{thm}

The number of minimal subsets is  related to the cardinality of sensitive sets. For example,
it was shown in \cite{Ye2008} that
if a transitive system has a dense set of minimal points but is not minimal, then
there exists an infinite sensitive set.
Moreover, if there are uncountable pairwise disjoint minimal subsets, then there exists an uncountable sensitive set.

In 2011, Huang, Lu and Ye got a fine structure of sensitive sets in minimal systems.

\begin{thm}[\cite{Huang2011}]
Let $(X, T)$ be a minimal dynamical system and $\pi:(X,T)\to (Y,S)$ be the factor map to
the maximal equicontinuous factor of $(X,T$). Then
\begin{enumerate}
  \item each sensitive set of $(X,T)$ is contained in some $\pi^{-1}(y)$ for some $y\in Y$;
  \item for each $y\in Y$, $\pi^{-1}(y)$  is a sensitive set of $(X,T$).
\end{enumerate}
Consequently, $(X, T)$ is $n$-sensitive, not $(n + 1)$-sensitive, if and only if
$\max_{y\in Y}\#(\pi^{-1}(y)) = n$.
\end{thm}

\subsection{Li-Yorke sensitivity}
A concept that combines sensitivity and Li-Yorke $\delta$-scrambled pairs was proposed by Akin and Kolyada in~\cite{Akin2003}, which is called Li-Yorke sensitivity.

\begin{defn}
A dynamical system $(X,T)$ is called \emph{Li-Yorke sensitive}
if there exists some $\delta>0$ such that for any $x\in X$ and $\ep>0$,
there is $y\in X$ satisfying $d(x,y)<\ep$ such that
\[\liminf_{n\to\infty} d(T^nx,T^n y)=0\  \text{and}\   \limsup_{n\to\infty} d(T^nx,T^n y)>\delta.\]
\end{defn}

Let $(X,T)$ be a dynamical system and $x\in X$.
The \emph{proximal cell} of $x$ is the set $\{y\in X:\colon (x,y) \emph{ is proximal}\,\}$.
First, it was proved in~\cite{Keynes1969} that for a weakly mixing $(X,T)$ the set of points $x$ at which
the proximal cell is residual in $X$ is itself residual in $X$.
Later it was proved by Furstenberg in~\cite{Furstenberg1981} that
the proximal cell of every point is residual, provided that the system is minimal and
weakly mixing.
Finally, Akin and Kolyada proved in~\cite{Akin2003} that
in any weakly mixing system the proximal cell of every point is residual.
Note that the authors in~\cite{Huang2004a} got more about the structure of the proximal cells
of $\F$-mixing systems, where $\F$ is an Furstenberg family.

In \cite{Akin2003}, Akin and Kolyada also proved that
\begin{thm}[\cite{Akin2003}]
If  a non-trivial dynamical system $(X,T)$ is   weakly mixing then it is Li-Yorke sensitive.
\end{thm}
See a survey paper~\cite{Kolyada2004} for more results about Li-Yorke sensitivity
and its relation to other concepts of chaos.
But the following question remains open.

\begin{prob}
Are all Li-Yorke sensitive systems Li-Yorke chaotic?
\end{prob}

\subsection{Mean equicontinuity and mean sensitivity}

\begin{defn}
A dynamical system $(X,T)$ is called \emph{mean equicontinuous} if for every $\ep>0$,
there exists some $\delta>0$ such that whenever $x,y\in X$ with $d(x,y)<\delta$,
\[\limsup_{n\to\infty}\frac{1}{n}\sum_{i=0}^{n-1}d(T^ix,T^iy)<\ep.\]
\end{defn}

This definition is equivalent to the notion of mean-L-stability which was first introduced by Fomin \cite{Fomin1951}.
It is an open question if every ergodic invariant measure on a mean-L-stable
system has discrete spectrum \cite{Scarpellini1982}.
The authors in \cite{Li2013b} firstly gave an affirmative answer to this question.
See \cite{Garcia-Ramos2014} for another approach to this question.

\begin{thm}[\cite{Li2013b}]
If a dynamical system $(X,T)$ is mean equicontinuous, then every ergodic invariant measure on $(X,T)$ has discrete spectrum
and hence the topological entropy of $(X,T)$ is zero.
\end{thm}

Similarly, we can define the local version and the opposite side of mean equicontinuity.

\begin{defn}
Let $(X,T)$ be a dynamical system.
A point $x\in X$ is called \emph{mean equicontinuous} if for every $\ep>0$,
there exists some $\delta>0$ such that for every $y\in X$ with $d(x,y)<\delta$,
\[\limsup_{n\to\infty}\frac{1}{n}\sum_{i=0}^{n-1}d(T^ix,T^iy)<\ep.\]
A transitive system is called \emph{almost mean equicontinuous} if there is at least one mean equicontinuous point.
\end{defn}

\begin{defn}
A dynamical system $(X,T)$ is called \emph{mean sensitive}
there exists some $\delta>0$ such that for every $x\in X$
and every neighborhood $U$ of $x$, there exists $y\in U$ and $n\in\N$ such that
\[\limsup_{n\to\infty}\frac{1}{n}\sum_{i=0}^{n-1}d(T^ix,T^iy)>\delta.\]
\end{defn}

We have the following dichotomy for transitive systems and minimal systems.

\begin{thm}[\cite{Li2013b}] \label{almostmeorms}
If a dynamical system $(X,T)$ is transitive, then $(X,T)$ is either almost mean equicontinuous or mean sensitive.
In particular, if $(X,T)$ is a minimal system, then $(X,T)$ is either mean equicontinuous or mean sensitive.
\end{thm}

Recall that an almost equicontinuous system is uniformly rigid and thus has zero topological
entropy. The following Theorem~\ref{p-entropy} shows that an almost mean equicontinuous
system behaves quite differently.

\begin{thm}[\cite{Li2013b}] \label{p-entropy}
In the full shift $(\Sigma_2,\sigma)$, every minimal subshift $(Y,\sigma)$ is contained
in an almost mean equicontinuous subshift $(X,\sigma)$.
\end{thm}

Since it is well-known that there are many minimal subshifts of $(\Sigma_2,\sigma)$ with positive
topological entropy, an immediate corollary of Theorem~\ref{p-entropy} is the following result.
\begin{cor}\label{cor:almost-mean-equi-entropy}
There exist many almost mean equicontinuous systems which have positive topological entropy.
\end{cor}

Globally speaking a mean equicontinuous system is `simple', since it is a Banach proximal
extension of an equicontinuous system and each of its ergodic measures has discrete
spectrum. Unfortunately, the local version does not behave so well, as Theorem~\ref{p-entropy}
shows. We will introduce the notion of Banach mean equicontinuity, whose
local version has the better behavior that we are looking for.

Let $(X,T)$ be a dynamical system. We say that $(X,T)$ is \emph{Banach mean equicontinuous}
if for every $\ep>0$,  there exists some $\delta>0$ such that whenever $x,y\in X$ with $d(x,y)<\delta$,
\[\limsup_{M-N\to\infty}\frac{1}{M-N}\sum_{i=N}^{M-1}d(T^ix,T^iy)<\ep.\]

A point $x\in X$ is called \emph{Banach mean equicontinuous}
if for every $\ep>0$, there exists some $\delta>0$ such that for every $y\in B(x,\delta)$,
\[\limsup_{M-N\to\infty}\frac{1}{M-N}\sum_{i=N}^{M-1}d(T^ix,T^iy)<\ep.\]
We say that a transitive system $(X,T)$ is \emph{almost Banach mean equicontinuous}
if there exists a transitive point which is Banach mean equicontinuous.

A dynamical system $(X,T)$ is \emph{Banach mean sensitive}
if there exists some $\delta>0$ such that for every $x\in X$ and every $\ep>0$
there is $y\in B(x,\ep)$ satisfying
\[\limsup_{M-N\to\infty}\frac{1}{M-N}\sum_{i=N}^{M-1}d(T^ix,T^iy)>\delta.\]

We also have the following dichotomy for transitive systems and minimal systems.

\begin{thm}[\cite{Li2013b}]\label{thm:BME-BMS}
If a dynamical system $(X,T)$ is transitive,
then $(X,T)$ is either almost Banach mean equicontinuous or Banach mean sensitive.
If $(X,T)$ is a minimal system, then $(X,T)$ is either  Banach mean equicontinuous or Banach mean sensitive.
\end{thm}

\begin{thm}[\cite{Li2013b}] \label{thm:BMS}
Let $(X,T)$ be a transitive system.
If the topological entropy of $(X,T)$ is positive, then $(X,T)$ is Banach mean sensitive.
\end{thm}

Combining Theorems~\ref{thm:BME-BMS} and~\ref{thm:BMS}, we have the following corollary.
\begin{cor}\label{cor:ABME}
If $(X,T)$ is almost Banach mean equicontinuous then the topological entropy
of $(X,T)$ is zero.
\end{cor}

By Corollary~\ref{cor:almost-mean-equi-entropy},
there are many almost mean equicontinuous systems which have
positive entropy.
Then by Corollary~\ref{cor:ABME} they are not almost Banach mean equicontinuous.
But the following question is still open.
\begin{prob}
Is there a minimal system which is mean equicontinuous but not Banach
mean equicontinuous?
\end{prob}

\subsection{Other generalization of sensitivity}

Note that sensitivity can be generalized in other ways,
see \cite{Moothathu2007, Wang2010, Huang2014b} for example.

Let $(X,T)$ be a dynamical system. For $\delta>0$ and a non-empty open subset $U$, let
\begin{align*}
  N(\delta,U)&=\{n\in\N\colon \exists x,y \in U \ \text{with} \ d(T^nx,T^ny)>\delta\}\\
  &= \{n\in\N\colon\diam(T^n(U))>\delta\}.
\end{align*}

Let $\F$ be a Furstenberg family.
According to Moothathu \cite{Moothathu2007}, we say that a dynamical system $(X,T)$ is \emph{$\F$-sensitive}
if there exists some $\delta>0$ such that for any non-empty open subset $U$ of $X$, $N(\delta,U)\in \F$.
$\F$-sensitivity for some special families were discussed
in~\cite{Tan2011a, Garcia-Ramos2014, Liu2014, LiR2014}.

A dynamical system $(X,T)$ is called \emph{multi-sensitive} if there
exists some $\delta>0$ such that for any finite open non-empty subsets $U_1,\ldots, U_n$ of $X$,
\[\bigcap_{i=1}^n N(\delta, U_i)\not=\emptyset.\]
In \cite{Huang2014b} among other things, Huang, Kolyada and
Zhang proved that for a minimal system thick sensitivity is equivalent to multi-sensitivity.
Moreover, they showed the following dichotomy for minimal systems.
\begin{thm}\label{20141122}
Let $(X, T)$ be a minimal system. Then $(X, T)$ is multi-sensitive
if and only if it is not an almost one-to-one extension of its maximal equicontinuous factor.
\end{thm}

Results similar to Theorem \ref{20141122} for other families will be appeared in~\cite{Ye2014} by Ye and Yu.

\section{Chaos in transitive systems}\label{section4}

In this section, we discuss various kinds of chaos in transitive systems.
\subsection{Auslander-Yorke chaos}
In \cite{Auslander1980} Auslander and Yorke defined a kind of chaos as
``topological transitivity plus pointwise instability''.
This leads to the following definition of Auslander-Yorke chaos.

\begin{defn}  %[\cite{Auslander1980}]
A dynamical system $(X,T)$ is called \emph{Auslander-Yorke chaotic} if it both transitive and sensitive.
\end{defn}

Due to Ruelle and Takens' work on turbulence~\cite{Ruelle1971},
Auslander-Yorke chaos was also called \emph{Ruelle-Takens chaos}
(see \cite{Xiong2005} for example).
It should be noticed that there is no implication relation between Li-Yorke chaos and Auslander-Yorke chaos.
For example, any non-equicontinuous distal minimal system is sensitive, so
it is  Auslander-Yorke chaotic, but it has  no Li-Yorke scrambled pairs.
On the other hand, there are non-periodic transitive systems with a fixed point
that are not sensitive~\cite{Akin2001}; by Theorem~\ref{thm:HY02} they are Li-Yorke chaotic.

The work of Wiggins~\cite{Wiggins1990} leads to the following definition.

\begin{defn}
A dynamical system $(X,T)$ is called \emph{Wiggins chaotic} if
there exists a subsystem $(Y,T)$ of $(X,T)$ such that
$(Y,T)$ is both transitive and sensitive.
\end{defn}

In~\cite{Ruette2005}, Ruette investigated transitive and sensitive subsystems for
interval maps. She showed that
\begin{thm}
\begin{enumerate}
\item For a transitive map $f\colon [0,1]\to[0,1]$, if it is  Wiggins chaotic then it is also Li-Yorke chaotic.
\item There exists a continuous map $f\colon [0, 1]\to [0, 1]$ of zero topological entropy which is Wiggins chaotic.
\item There exists a Li-Yorke chaotic continuous map $f\colon[0, 1]\to [0, 1]$ which is not Wiggins chaotic.
\end{enumerate}
\end{thm}

\subsection{Devaney chaos}
In his book~\cite{Devaney1989}, Devaney  proposed a new kind of chaos,
which is usually called Devaney chaos.
\begin{defn}
A dynamical system $(X,T)$ is called \emph{Devaney chaotic} if it satisfies the following three properties:
\begin{enumerate}
  \item $(X,T)$ is transitive;
  \item $(X,T)$ has sensitive dependence on initial conditions;
  \item the set of periodic points of $(X,T)$ is dense in $X$.
\end{enumerate}
\end{defn}

Sensitive dependence is widely understood as the
central idea in Devaney chaos, but it is implied by transitivity and
density of periodic points, see~\cite{Banks1992} or~\cite{Glasner1993}.
In 1993, S.~Li~\cite{Li1993} showed that

\begin{thm}
Let $f:[0,1]\to[0,1]$ be a continuous map.
Then $f$ has positive topological entropy if and only if there exists a Devaney chaotic subsystem.
\end{thm}

In 1996, Akin, Auslander and Berg~\cite{Akin1996} discussed in details when a transitive map is sensitive.
An important question is that: does Devaney chaos imply Li-Yorke chaos?
In 2002, Huang and Ye gave a positive answer to this  long-standing open problem.

\begin{thm}[\cite{Huang2002}]\label{thm:HY02}
Let $(X,T)$ be a non-periodic transitive system.
If there exists a periodic point, then it is Li-Yorke chaos.
Particularly,  Devaney chaos implies Li-Yorke chaos.
\end{thm}
A key result in the proof of Theorem~\ref{thm:HY02} is called \emph{Huang-Ye equivalences}
in~\cite{Akin2003}.
To state this result, we need some notations.

For every $\ep>0$, let $\overline{V_\ep}=\{(x,y)\in X\times X\colon d(x,y)\leq\ep \}$ and
\[Asym_\ep(T)= \bigcup_{n=0}^\infty\bigg(\bigcap_{k=n}^\infty T^{-k} \overline{V_\ep}\bigg).\]
For every $x\in X$,  let
\[Asym_\ep(T)(x)=\big\{y\in X\colon (x,y)\in Asym_\ep(T)\big\}.\]
\begin{thm}[Huang-Ye Equivalences]
For a dynamical system $(X, T)$ the following conditions are equivalent.
\begin{enumerate}
\item $(X,T)$ is sensitive.
\item There exists $\ep>0$ such that $Asym_\ep(T)$ is a first category subset of $X\times X$.
\item There exists $\ep>0$ such that for every $x\in X$, $Asym_\ep(T)(x)$ is
a first category subset of $X$.
\item There exists $\ep>0$ such that for every $x\in X$, $x\in\overline{X\setminus Asym_\ep(T)(x)}$.
\item There exists $\ep>0$ such that the set
\[\big\{(x,y)\in X\times X\colon \limsup_{n\to\infty} d(T^nx,T^ny)>\ep\big\}\]
is dense in  $X\times X$.
\end{enumerate}
\end{thm}

Note that the first part of Theorem~\ref{thm:HY02} does not involve the sensitive property directly.
In 2004, Mai proved that Devaney chaos implies the existence of $\delta$-scrambled sets for some $\delta>0$.
Using a method of direct construction, he showed that
\begin{thm}[\cite{Mai2004}]
Let $(X,T)$ be a non-periodic transitive system.
If there exists a periodic point of period $p$, then
there exist Cantor sets $C_1\subset C_2\subset \dotsb$ in $X$ consisting of
transitive points such that
\begin{enumerate}
\item each $C_n$ is a synchronously proximal set, that is $\liminf\limits_{k\to\infty} \diam(T^k C_n)=0$;
\item $S:=\bigcup_n C_n$ is scrambled, and $\bigcup_{i=0}^{p-1}T^iS$ is dense in $X$;
\item furthermore, if $(X,T)$ is sensitive, then $S$ is $\delta$-scrambled  for some $\delta>0$.
\end{enumerate}
\end{thm}

\subsection{Multivariant Li-Yorke chaos} \label{sec:multi-LY-chaos}
The scrambled set in Li-Yorke chaos only compares the trajectories of two points,
It is natural to consider the trajectories of finite points.
In 2005, Xiong introduced the following multivariant chaos in the sense of Li-Yorke.

\begin{defn}[\cite{Xiong2005}]
Let $(X,T)$ be a dynamical system and $n\geq 2$.
A tuple $(x_1,x_2,\ldots,\allowbreak x_n)\in X^n$ is called \emph{$n$-scrambled} if
\[\liminf_{k\to\infty}\max_{1\leq i<j\leq n}d\big(f^k(x_i), f^k(x_j)\big)=0\]
and
\[\qquad\limsup_{k\to\infty}\min_{1\leq i<j\leq n}d\big(f^k(x_i), f^k(x_j)\big)>\delta>0.\]
A subset $C$ of $X$ is called \emph{$n$-scrambled} if any pairwise distinct $n$ points in $C$ form an $n$-scrambled tuple.
The dynamical system $(X,T)$ is called \emph{Li-Yorke $n$-chaotic} if
there exists an uncountable $n$-scrambled set.

Similarly, if the separated constant $\delta$ is uniform for all pairwise distinct $n$-tuples in $C$,
we can define \emph{$n$-$\delta$-scrambled sets} and \emph{Li-Yorke $n$-$\delta$-chaos}.
\end{defn}

\begin{thm}[\cite{Xiong2005}]\label{thm:X2005}
Let $(X,T)$ be a non-trivial transitive system.
If $(X,T)$ has a fixed point, then there exists a dense Mycielski subset $C$ of $X$ such that
$C$ is $n$-scrambled for all $n\geq 2$.
\end{thm}

%\begin{rem}\label{rem:X2005-n-sensitive}
%In Theorem~\ref{thm:X2005}, if in addition $(X,T)$ is $m$-sensitive,
%then the set $C$ can be chosen to be Li-Yorke $m$-$\delta_m$-scrambled
%with some $\delta_m>0$ \cite{Xiong2005}.
%\end{rem}
%
%\begin{thm}[\cite{Ye2008}]\label{thm:YZ2008}
%If $(X,T)$ is Devaney chaotic, then it is $n$-sensitive for any $n\geq 2$.
%\end{thm}
%
%Combining Theorems~\ref{thm:X2005}, \ref{thm:YZ2008} and Remark~\ref{rem:X2005-n-sensitive},
%we have the following result.
%\begin{thm}
%Let $(X,T)$ be a dynamical system.
%If $(X,T)$ is Devaney chaotic, then there exists a Mycielski subset $C$ of $X$ such that for any $n\geq 2$,
%$C$ is $n$-$\delta_n$-scrambled with some $\delta_n>0$.
%Moreover, if $(X,T)$ has a fixed point, then $C$ can be chosen to be dense in $X$.
%\end{thm}

In 2011, Li proved that there is no $3$-scrambled tuples for an interval map with zero topological entropy.

\begin{thm}[\cite{Li2011}]
Let $f:[0,1]\to[0,1]$ be a continuous map.
If $f$ has zero topological entropy, then there is no  $3$-scrambled tuples.
\end{thm}

Note that there exists a continuous map $f\colon[0,1]\to [0,1]$,
which is Li-Yorke chaotic but has zero topological entropy (see Theorem~\ref{thm:LY-chaos-0entropy}).
Then we have the following result, which shows that Li-Yorke $2$-chaos and $3$-chaos are essentially different.

\begin{cor}
There exists a Li-Yorke $2$-chaotic system which has no $3$-scrambled tuples.
\end{cor}

\subsection{Uniform chaos}
Recall that a dynamical system $(X,T)$ is \emph{scattering} if for any minimal system $(Y,S)$
the product system $(X\times Y,T\times S)$ is transitive.
It is not hard to show that every weakly mixing system is  scattering.
In~\cite{Huang2002}, Huang and Ye also proved that scattering implies Li-Yorke chaos.

\begin{thm}%[\cite{Huang2002}]
If a non-trivial  dynamical system $(X,T)$ is scattering, then there is a dense Mycielski scrambled set.
\end{thm}

Since every system has a minimal subsystem,
for a scattering system $(X,T)$ there exists a subsystem $(Y,S)$ of $(X,T)$
such that $(X\times Y,T\times T)$ is transitive. In 2010, Akin \textit{et al.}
showed that a dynamical system satisfying this property is
more complicated than Li-Yorke chaos, and proposed the concept of uniform chaos~\cite{Akin2010}.

\begin{defn}
Let $(X,T)$ be a dynamical system. A subset $K$ of $X$ is said to be
\begin{enumerate}[(1)]
\item \emph{uniformly recurrent} if for every $\ep>0$ there exists an $n\in\N$
with  $d(T^nx,x)<\ep$ for all $x\in K$;
\item \emph{recurrent} if every finite subset of $K$ is uniformly recurrent;
\item \emph{uniformly proximal} if for every $\ep>0$ there exists an $n\in\N$
with $\diam(T^n K)<\ep$;
\item \emph{proximal} if  every finite subset of $K$ is uniformly proximal.
\end{enumerate}
\end{defn}

\begin{defn}
Let $(X,T)$ be a dynamical system. A subset $K\subset X$ is called
a \emph{uniformly chaotic set} if there are Cantor sets $C_1\subset C_2\subset \dotsb$ such that
\begin{enumerate}
\item for each $N\in\N$, $C_N$ is uniformly recurrent;
\item for each $N\in\N$, $C_N$ is uniformly proximal;
\item $K:=\bigcup_{i=1}^\infty C_n$ is a recurrent subset of $X$ and also a proximal subset of $X$.
\end{enumerate}
The system $(X,T)$ is called (\emph{densely}) \emph{uniformly chaotic}
if it has a (dense) uniformly chaotic subset of $X$.
\end{defn}
By the definition, any uniformly chaotic set is scrambled, then uniform chaos implies Li-Yorke chaos.
The following is the main result of~\cite{Akin2010}.
\begin{thm}%[\cite{Akin2010}]
Let $(X,T)$ be a non-trivial dynamical system.
If there exists a subsystem $(Y,T)$ of $(X,T)$ such that  $(X\times Y, T\times T)$ is transitive,
then $(X,T)$ is densely uniformly chaotic.
\end{thm}

As a corollary, we can show that many transitive systems are uniformly chaotic.
Note that
a dynamical system is \emph{weakly scattering} if its product with any minimal equicontinuous system is
transitive.

\begin{cor}\label{cor:unf-chaos}
If $(X,T)$ is a dynamical system without isolated points and
satisfies one of the following properties,
then it is densely uniformly chaotic:
\begin{enumerate}
\item $(X,T)$ is transitive and has a fixed point;
\item $(X,T)$ is totally transitive with a periodic point;
\item $(X,T)$ is scattering;
\item $(X,T)$ is weakly scattering with an equicontinuous minimal subsystem;
\item $(X,T)$ is weakly mixing.
\end{enumerate}
Furthermore, if $(X,T)$ is transitive and has a periodic point of period $p$,
then there is a closed $T^p$-invariant subset $X_0$ of $X$,
such that $(X_0,T^p)$ is densely uniformly chaotic and $X=\bigcup_{j=0}^{p-1}T^j X_0$;
In particular, $(X,T)$ is uniformly chaotic.
\end{cor}

\begin{rem}
By Corollary~\ref{cor:unf-chaos}, Devaney chaos implies uniform chaos.
\end{rem}

\section{Various types of distributional chaos and its generalization}\label{section5}

\subsection{Distributional chaos}
In~\cite{Schweizer1994}, Schweizer and Sm{\'{\i}}tal used ideas from probability theory to
develop a new definition of chaos, so called \emph{distributional chaos}.
Let $(X,T)$ be a dynamical system.
For a pair of points $x,y$ in $X$ and a positive integer $n$, we define a function
$\Phi^n_{xy}$ on the real line by
\[ \Phi^n_{xy}(t)=\tfrac 1 n \#\big\{0\leq i \leq n-1\colon d(T^ix, T^iy)<t\big\},\]
where $\#(\cdot)$ denotes the number of elements of a set.
Clearly, the function $\Phi^n_{xy}$ is non-decreasing, has minimal value $0$
(since $\Phi^n_{xy}(t)=0$ for all $t\leq 0$), has maximum value $1$
(since $\Phi^n_{xy}(t)=1$ for all $t$ greater then the diameter of $X$),
and is left continuous. Then $\Phi^n_{xy}$ is a distribution function
whose value at $t$ may be interpreted as the probability that the distance between
the initial segment of length $n$ of the trajectories of $x$ and $y$ is less than $t$.

We are interested in the asymptotic behavior of the function $\Phi^n_{xy}$
as $n$  gets large. To this end,  we consider the functions
$\Phi^*_{xy}$ and $\Phi_{xy}$ defined by
\[\Phi^*_{xy}(t)=\limsup_{n\to\infty} \Phi^n_{xy}(t)\quad\text{and}\quad
\Phi_{xy}(t)=\liminf_{n\to\infty} \Phi^n_{xy}(t).\]
Then $\Phi^*_{xy}$ and $\Phi_{xy}$ are distribution functions with
$\Phi_{xy}(t)\leq \Phi_{xy}^*(t)$ for all $t$.
It follows that $\Phi^*_{xy}$ is an asymptotic measure of how close $x$ and $y$ can come together,
while $\Phi_{xy}$  is an asymptotic measure of their maximum separation.
We shall refer to $\Phi^*_{xy}$ as the \emph{upper distribution},
and to $\Phi_{xy}$ as the \emph{lower distribution} of $x$ and $y$.

\begin{defn}
Let $(X,T)$ be a dynamical system.
A pair $(x,y)\in X\times X$ is called \emph{distributionally scrambled} if
\begin{enumerate}
\item for every $t>0$, $\Phi^*_{xy}(t)=1$, and
\item there exists some $\delta>0$ such that  $\Phi_{xy}(\delta)=0$.
\end{enumerate}
A subset $C$ of $X$ is called \emph{distributionally scrambled}
if any two distinct points in $C$ form a distributionally scrambled pair.
The dynamical system $(X,T)$ is called \emph{distributionally chaotic} if there
exists an uncountable distributionally scrambled set.

Similarly, if the separated constant $\delta$ is uniform for all non-diagonal pairs in $C$,
we can define \emph{distributionally $\delta$-scrambled sets} and \emph{distributional $\delta$-chaos}.
\end{defn}

According to the definitions, it is clear that any distributionally scrambled pair is scrambled,
and then  distributional chaos is stronger than Li-Yorke chaos.
In~\cite{Schweizer1994}, Schweizer and Sm{\'{\i}}tal showed that

\begin{thm}[\cite{Schweizer1994}]
Let $f\colon[0,1]\to [0,1]$ be a continuous map. Then $f$ is distributionally chaotic if and only if
it has positive topological entropy.
\end{thm}

In 1998, Liao and Fan~\cite{Liao1998} constructed a minimal system with zero topological entropy which is
distributionally chaotic.
In 2006, Oprocha~\cite{Oprocha2006} showed that weak mixing and Devaney chaos
do not imply distributional chaos.

In Section~\ref{sec:multi-LY-chaos}, Li-Yorke chaos is extended to a multivariant version.
In fact, it is clear that any kind of chaos defined by scrambled pairs
can be extended to multivariant version.

In \cite{Tan2014}, Tan and Fu showed that distributionally $n$-chaos and $(n+1)$-chaos are essentially different.
\begin{thm}[\cite{Tan2014}]
For every $n\geq 2$, there exists a transitive system which is distributionally $n$-chaotic
but without any distributionally $(n+1)$-scrambled tuples.
\end{thm}

In 2013, Li and Oprocha showed that
\begin{thm}[\cite{Li2013O}] \label{thm:n-chaos-example}
For every $n\geq 2$, there exists a dynamical system which is distributionally $n$-chaotic
but not Li-Yorke $(n+1)$-chaotic.
\end{thm}

Note that the example constructed in the proof of Theorem~\ref{thm:n-chaos-example}
contains $(n+1)$-scrambled tuples. Then, a natural question is as following.
\begin{prob} \label{prob:DC1-3scrambled-tuples}
Is there a dynamical system $(X,T)$ with an uncountable distributionally $2$-scrambled
set but without any $3$-scrambled tuples?
\end{prob}

Recently, Dole{\v{z}}elov{\'a} made some progress on the Problem \ref{prob:DC1-3scrambled-tuples},
but the original problem is still open. She showed that
\begin{thm}[\cite{Dolezelova2014}]
There exists a dynamical system $X$ with an infinite extremal distributionally scrambled set
but without any scrambled triple.
\end{thm}
\begin{thm}[\cite{Dolezelova2014}]
There exists an invariant Mycielski (not closed) set $X$ in the full shift with
an uncountable extremal distributionally $2$-scrambled set but without any $3$-scrambled tuple.
\end{thm}

\subsection{The three versions of distributional chaos}
Presently, we have at least three different definitions of
distributionally scrambled pair (see~\cite{Smital2004} and~\cite{Balibrea2005}).
The original distributionally scrambled pair is said to be \emph{distributionally scrambled of type 1}.

\begin{defn}
A pair $(x,y)\in X\times X$ is called \emph{distributionally scrambled of type 2} if
\begin{enumerate}
  \item for any $t>0$, $\Phi^*_{xy}(t)=1$, and
  \item there exists some $\delta>0$ such that $\Phi_{xy}(\delta)<1$.
\end{enumerate}

A pair $(x,y)\in X\times X$ is called \emph{distributionally scrambled of type 3} if
there exists an interval $[a,b]\subset (0,\infty)$ such that $ \Phi_{xy}(t)<\Phi^*_{xy}(t)$ for all $t\in [a,b]$.

A subset $C$ of $X$ is called \emph{distributionally scrambled of type i} ($i=1,2,3$) if
any two distinct points in $C$ form a distributionally scrambled pair of type $i$.
The dynamical system $(X,T)$ is called \emph{distributionally chaotic of type i} (DC{$i$} for short)
if there exists an uncountable distributionally scrambled set of type $i$.
\end{defn}

It should be noticed that in~\cite{Smital2004} and~\cite{Balibrea2005}
DC$i$ chaos only requires the existence of one distributionally scrambled pair of type $i$.
It is not hard to see that any distributionally scrambled pair of type 1 or 2
is topological conjugacy invariant~\cite{Smital2004}.
But a distributionally scrambled pair of type 3 may not be topological conjugacy invariant~\cite{Balibrea2005}.
It is also not hard to construct a dynamical system has a distributionally scrambled pair of type 2,
but no distributionally scrambled pairs of type 1.
A really interesting example constructed in~\cite{Tan2009} shows that
there exists a minimal subshift which is distributionally chaotic  of type 2, while it
does not contain any distributionally scrambled pair of type 1.

\subsection{The relation between distributional chaos and positive topological entropy}
It is known that in general there is no implication relation between DC1 and
positive topological entropy (see~\cite{Liao1998} and~\cite{Pikula2007}).
In \cite{Smital2006}~Sm\'{\i}tal conjectured that positive topological entropy does imply DC2.
Oprocha showed that this conjecture is true for uniformly positive entropy minimal systems~\cite{Oprocha2011a}.
Finally, Downarowicz solved this problem by proving this conjecture in general case~\cite{Downarowicz2011}.

\begin{thm}[\cite{Downarowicz2011}]\label{thm:D11}
If a dynamical system $(X,T)$ has positive topological entropy,
then there exists a Cantor distributional $\delta$-scrambled set of type $2$ for some $\delta>0$.
\end{thm}

It was observed in~\cite{Downarowicz2011} that a pair $(x,y)$ is DC2-scrambled if and only if
it is mean scrambled in the Li-Yorke sense, that is
\[\liminf_{N\to \infty}\ \frac{1}{N}\sum_{k=1}^Nd\big(T^k x , T^k y)=0\]
and
\[\limsup_{N\to \infty}\ \frac{1}{N}\sum_{k=1}^Nd\big(T^k x , T^k y)>0.\]
Recently, Huang, Li and Ye showed that positive topological entropy implies a multivariant version of mean Li-Yorke chaos.
\begin{thm}[\cite{Huang2014}]\label{thm:PTE-n-MLC}
If a dynamical system $(X,T)$ has positive topological  entropy,
then there exists a Mycielski multivariant mean Li-Yorke $(\delta_n)$-scrambled subset $K$ of $X$, that is
for every $n\geq 2$ and every $n$ pairwise distinct points  $x_1,\dotsc,x_n\in K$, we have
\[\liminf_{N\to \infty}\ \frac{1}{N}\sum_{k=1}^N\max_{1\leq i<j\leq n}d\big(T^k x_i , T^k x_j)=0\ \]
and
\[\limsup_{N\to \infty}\ \frac{1}{N}\sum_{k=1}^N\min_{1\leq i<j\leq n}d\big(T^k x_i , T^k x_j)>\delta_n>0.\]
\end{thm}

We remark that the key tool in the proof of Theorem~\ref{thm:PTE-n-MLC}
is the excellent partition constructed in \cite[Lemma 4]{Blanchard2002a},
which is different from the one used in~\cite{Downarowicz2011}.
So among other things for $n=2$ the authors also obtained a new proof of Theorem~\ref{thm:D11}.

Note that in~\cite{Blanchard2002a} Blanchard, Host and Ruette showed that
the closure of the set of asymptotic pairs contains the set of entropy pairs.
Kerr and Li ~\cite{Kerr2007} proved that the intersection of the set of scrambled pairs
with the set of entropy pairs is dense in the set of entropy pairs.
We extended the above mentioned results to the following result.

\begin{thm}[\cite{Huang2014}]
If a dynamical system  has positive topological  entropy, then for any $n\geq 2$,
\begin{enumerate}
  \item the intersection of the set of asymptotic $n$-tuples with the set of topological entropy $n$-tuples
  is dense in the set of topological entropy $n$-tuples;
  \item the intersection of the set of mean $n$-scrambled tuples with the set of topological entropy $n$-tuples
  is dense in the set of topological entropy $n$-tuples.
\end{enumerate}
\end{thm}
\subsection{Chaos via Furstenberg families}
In Section~\ref{sec:pre}, we have shown that there is a powerful connection between
topological dynamics and Furstenberg families.
In 2007, Xiong, L\"u and Tan introduced the notion of chaos via  Furstenberg families.
It turned out that the Li-Yorke chaos and some version of distributional chaos can
be described as chaos in Furstenberg families sense.

\begin{defn}[\cite{Xiong2007, Tan2009}]
Let $(X,T)$ be a dynamical system. Let
$\F_1$ and $\F_2$ be two Furstenberg families.
A pair $(x,y)\in X\times X$ is called \emph{$(\F_1,\F_2)$-scrambled} if it satisfies
\begin{enumerate}
\item for any $t>0$, $\{n\in\Z\colon d(T^nx,T^ny)<t\}\in \F_1$, and
\item there exists some $\delta>0$ such that $\{n\in\Z\colon d(T^nx,T^ny)>\delta\}\in \F_2$.
\end{enumerate}
A subset $C$ of $X$ is called \emph{$(\F_1,\F_2)$-scrambled} if any two distinct
points in $C$ form a $(\F_1,\F_2)$-scrambled pair.
The dynamical system $(X,T)$ is called \emph{$(\F_1,\F_2)$-chaotic} if
there is an uncountable $(\F_1,\F_2)$-scrambled set in $X$.

Similarly, if the separated constant $\delta$ is uniform for all non-diagonal pairs in $C$,
we can define \emph{$(\F_1,\F_2)$-$\delta$-scrambled sets} and \emph{$(\F_1,\F_2)$-$\delta$-chaos}.
\end{defn}

Let $F$ be a subset of $\Z$. The \emph{upper density} of $F$ is defined as
\[\overline{d}(F)=\limsup_{n\to\infty}\frac 1 n \#(F\cap\{0,1,\dots,n-1\}).\]
For any $a\in [0,1]$, let
\[\M(a)=\{F\subset\Z\colon F\text{ is infinite and }\overline{d}(F)\geq a\}.\]
Clearly, every $\M(a)$ is a Furstenberg family.
The following proposition can be easily verified by the definitions,
which shows that Li-Yorke chaos and distributional chaos can be characterized via
Furstenberg families.

\begin{prop}
Let $(X,T)$ be a dynamical system and $(x,y)\in X\times X$. Then
\begin{enumerate}[(1)]
  \item $(x,y)$ is scrambled if and only if it is $(\M(0),\M(0))$-scrambled.
  \item $(x,y)$ is distributionally scrambled of type 1 if and only if it is $(\M(1),\M(1))$-scrambled.
  \item $(x,y)$ is distributionally scrambled of type 2 if and only if it is
  $(\M(1),\M(t))$-scrambled with some $t>0$.
\end{enumerate}
\end{prop}

A Furstenberg family $\F$ is said to be \emph{compatible} with the dynamical system $(X,T)$,
if for every non-empty open subset $U$ of $X$
the set $\{x\in X\colon\{n\in\Z\colon T^n x\in U\}\in\F\,\}$ is a $G_\delta$ subset of $X$.
\begin{lem}[\cite{Xiong2007}]
If $\F=\M(t)$ for some $t\in [0,1]$, then $\F$
is compatible with any dynamical system $(X,T)$.
\end{lem}

We have the following criteria for chaos via compatible Furstenberg families.
\begin{thm}[\cite{Xiong2007,Tan2009}]
Let $(X,T)$ be a dynamical system.
Suppose that there exists a fixed point $p\in X$ such that the set $\bigcup_{i=1}^\infty T^{-i}(p)$
is dense in $X$, and a non-empty closed invariant subset $A$ of $X$ disjoint with the point $p$
such that $\bigcup_{i=1}^\infty T^{-i}A$ is dense in $X$.
Then for every two Furstenberg families $\F_1$ and $\F_2$ compatible with the system $(X\times X,T\times T)$,
there exists a dense Mycielski $(\F_1,\F_2)$-$\delta$-scrambled set for some $\delta>0$.
In particular, the dynamical system $(X, T)$ is $(\M(1),\M(1))$-$\delta$-chaotic for some $\delta>0$.
\end{thm}

There are two important Furstenberg families: the collections of all syndetic sets and all sets with Banach density one.
They are so specifical that we use terminologies: syndetically scrambled sets and
Banach scrambled sets, when studying chaos via those two Furstenberg families.

Recall that a subset $F$ of $\Z$ is \emph{syndetic} if
there is $N\in\N$ such that $[n,n + N] \cap F\neq\emptyset$ for every  $n \in\Z$.
Let $(X,T)$ be a dynamical system.
A pair of points $(x,y)\in X^2$ is called \emph{syndetically proximal}
if for every $\eps>0$ the set $\{n\in\Z\colon d(T^nx,T^ny)<\eps\}$ is syndetic.
If $(X,T)$ is proximal, that is every pair in $X^2$ is proximal,
then every pair in $X^2$ is also syndetically proximal~\cite{Moothathu2011}.
Moothathu studied the syndetically proximal relation,
and identified certain sufficient conditions
for the syndetically proximal cell of each point to be small~\cite{Moothathu2011}.
A subset $S$ of $X$ is called \emph{syndetically scrambled}
if for any two distinct points $x,y\in S$,
$(x,y)$ is syndetically proximal but not asymptotic.
In~\cite{Moothathu2013} Moothathu and Oprocha systematically studied the syndetically proximal
relation and the possible existence of syndetically scrambled sets for many kinds of dynamical systems,
including various classes of transitive subshifts, interval maps, and topologically Anosov maps.

A subset $F$ of $\Z$ is said to have \emph{Banach density one}
if for every $\lambda < 1$ there exists $N \geq 1$ such that
$\#(F \cap I) \geq \lambda\#(I)$ for every subinterval $I$ of $\Z$ with $\#(I)\geq N$,
where $\#(I)$ denotes the number of elements of $I$.
A pair of points $(x, y) \in X^2$ is called \emph{Banach proximal} if for every $\eps > 0$
the set $\{n \in\Z\colon d(T^nx,T^ny)<\eps\}$ has Banach density one.
Clearly, every set with Banach density one is syndetic,
then a Banach proximal pair is syndetically proximal.
In~\cite{Li2014a} Li and Tu studied the structure of the Banach proximal relation.
Particularly, they showed that for a non-trivial minimal system
the Banach proximal cell of every point is of first category.

A subset $S$ of $X$  is called \emph{Banach scrambled}
if for any two distinct points $x,y\in S$,
$(x,y)$ is  Banach proximal but not asymptotic.
Note that it was shown in~\cite{Li2011} that for an interval map with zero topological
entropy, every proximal pair is Banach proximal. Moreover we have
\begin{thm}[\cite{Li2014a}]
Let $f\colon[0,1]\to[0,1]$ be a continuous map.
If $f$ is Li-Yorke chaotic, then it has a Cantor Banach scrambled set.
\end{thm}

\begin{thm}[\cite{Li2014a}]
Let $f\colon[0,1]\to[0,1]$ be a continuous map.
Then $f$ has positive topological entropy if and only if
there is Cantor set $S\subset [0,1]$ such that
$S$ is a Banach scrambled set and $f^n(S)\subset S$ for some $n\geq 1$.
\end{thm}

There exists a dynamical system with the whole space being a scrambled set
(see Theorems~\ref{thm:completely-scrambled-systems} and~\ref{thm:wm-uniformly-rigid}),
and then every pair is syndetically proximal.
It is shown in~\cite{Li2014a} that the whole space also can be a Banach scrambled set.

\begin{thm}[\cite{Li2014a}]
There exists a transitive system with the whole space being a Banach scrambled set.
\end{thm}

\section{Local aspects of mixing properties}

As we said before for a dynamical system with positive topological entropy, it is possible that there is no
weakly mixing subsystem. The notion of weakly mixing sets captures some complicated subset of the system.

\subsection{Weakly mixing sets}
In 1991, Xiong and Yang~\cite{Xiong1991} showed that for in a weakly mixing system
there are considerably many points in the domain whose orbits display
highly erratic time dependence.
More specifically, they showed that
\begin{thm}\label{thm:XY1991}
Let $(X,T)$ be a non-trivial dynamical system. Then
\begin{enumerate}
  \item $(X,T)$ is weakly mixing if and only if there exists a $c$-dense $F_\sigma$ subset $C$ of $X$
  satisfying for any subset $D$ of $C$ and any continuous map $f\colon D\to X$, there exists an
  increasing sequence $\{q_i\}$ of positive integers such that
  $\lim\limits_{i\to\infty} T^{q_i}x=f(x)$ for all $x\in D$;
  \item $(X,T)$ is strongly mixing if and only if for any increasing sequence $\{p_i\}$ of
  positive integers there exists a $c$-dense $F_\sigma$ subset $C$ of $X$
  satisfying for any subset $D$ of $C$ and any continuous map $f\colon D\to X$, there exists a
  subsequence $\{p_{n_i}\}$ of $\{p_i\}$ such that
  $\lim\limits_{i\to\infty} T^{p_{n_i}}x=f(x)$ for all $x\in D$.
\end{enumerate}
In particular, if $(X,T)$ is weakly mixing, then there exists some $\delta>0$
such that $(X,T)$ is Li-Yorke $\delta$-chaotic.
\end{thm}
Note that a subset $C$ is called \emph{$c$-dense} in $X$ if for every non-empty open subset
$U$ of $X$, the intersection $C\cap U$ has the cardinality of the continuum $c$.

Inspired by Xiong-Yang's result (Theorem~\ref{thm:XY1991}),
in 2008 Blanchard and Huang~\cite{Blanchard2008} introduced the concept of weakly mixing sets,
which can be regraded as a local version of weak mixing.

\begin{defn}
Let $(X,T)$ be a dynamical system. A closed  subset $A$ of $X$ is said to be \emph{weakly mixing}
if there exists a dense Mycielski subset $C$ of $A$ such that
 for any subset $D$ of $C$ and any continuous map $f\colon D\to A$, there exists an
increasing sequence $\{q_i\}$ of positive integers satisfying
$\lim\limits_{i\to\infty} T^{q_i}x=f(x)$ for all $x\in D$.
\end{defn}

Blanchard and Huang showed that in any positive entropy system there are many weakly mixing sets.
More precisely, let $WM(X, T)$ be the set of weakly mixing sets of $(X,T)$ and
$H(X, T)$ be the closure of the set entropy sets in the hyperspace.

\begin{thm}[\cite{Blanchard2008}]
If a dynamical system $(X,T)$ has positive topological entropy,
then the set $H(X, T )\cap WM(X, T )$ is a dense $G_\delta$ subset of $H(X,T)$.
\end{thm}

Moreover in \cite{Huang2008}, Huang showed that in any positive entropy system
there is a measure-theoretically ¡°rather big¡± set such that the closure of
the stable set of any point from the set contains a weakly mixing set.
Recall that the \emph{stable set} of a point $x\in X$ for $T$ is
\[W^s(X,T)=\{y\in X\colon \liminf_{n\to\infty}d(T^nx,T^ny)=0\}.\]

\begin{thm}[\cite{Huang2008}]
Let $(X,T)$ be a dynamical system and $\mu$ be an ergodic invariant measure on $(X,T)$ with positive entropy.
Then for $\mu$-a.e.\ $x\in X$ the closure of the stable set $\overline{W^s(X,T)}$ contains a weakly mixing set.
\end{thm}

Recently in \cite{Huang2014c}, Huang, Xu and Yi showed the existence of certain chaotic sets in the stable set
of positive entropy $G$-systems for certain countable, discrete, infinite
left-orderable amenable groups $G$.
We restate \cite[Theorem 1.2]{Huang2014c} in our setting as follows.

\begin{thm}[\cite{Huang2014c}]
Let $T$ be a homeomorphism on a compact metric space $X$.
If $\mu$ is an ergodic invariant measure on $(X,T)$ with positive entropy,
then there exists $\delta>0$ such that for $\mu$-a.e.\ $x\in X$
the stable set $W^s(X,T)$ contains a Cantor $\delta$-scrambled set for $T^{-1}$.
\end{thm}

\begin{rem}
The authors in~\cite{Blanchard2008} also discussed the relation between weakly mixing sets
and other chaotic properties.
\begin{enumerate}
\item Positive entropy is strictly stronger than the existence of weakly mixing sets,
which in turn is strictly stronger than Li-Yorke chaos.
\item There exists a Devaney chaotic system without any weakly mixing sets.
\end{enumerate}
\end{rem}

The following result is the well known Furstenberg intersection lemma, which shows that
weak mixing implies weak mixing of all finite orders.
\begin{lem}[\cite{Furstenberg1967}] \label{lem:Furstenberg-intersection-lemma}
If a dynamical system $(X,T)$ is weakly mixing,
then for any $k\in\N$ and any non-empty open subsets  $U_1,\dotsc,U_k$, $V_1,\dotsc,V_k$ of $X$,
\[ \bigcap_{i=1}^k N(U_i,V_i)\neq\emptyset.\]
\end{lem}

Using the idea in Lemma~\ref{lem:Furstenberg-intersection-lemma}
we can give the following alternative definition of weakly mixing sets.
\begin{prop}[\cite{Blanchard2008}]
Let $(X,T)$ be a dynamical system and let $A$ be a closed subset of $X$ but not a singleton.
Then  $A$ is  weakly mixing  if  and only if for any $k\in\N$ and any choice of non-empty
open subsets $U_1,\dotsc,U_k$, $V_1,\dotsc,V_k$
of $X$ intersecting $A$ (that is $A\cap U_i\neq\emptyset$, $A\cap V_i\neq\emptyset$ for $i=1,\dotsc,k$),
one has
\[ \bigcap_{i=1}^k N(U_i\cap A,V_i)\neq\emptyset.\]
\end{prop}

\begin{defn}
Let $(X,T)$ be a dynamical system, $A$ be a closed subset of $X$ but not a singleton and $k\geq 2$.
The set $A$ is said to be \emph{weakly mixing of order $k$} if for any choice of non-empty
open subsets $U_1,\dotsc,U_k$, $V_1,\dotsc,V_k$
of $X$  intersecting $A$, one has
\[ \bigcap_{i=1}^k N(U_i\cap A,V_i)\neq\emptyset.\]
Then $A$ is weakly mixing if and only if it is weakly mixing of order $k$ for all $k\geq 2$.
\end{defn}

In 2011, Oprocha and Zhang~\cite{Oprocha2011} studied weakly mixing sets of finite order,
and constructed an example showing that the concepts of weakly mixing sets of
order $2$ and of order $3$ are different.
They generalized this result to general weakly mixing sets of finite order in~\cite{Oprocha2013b}.

\begin{thm} [\cite{Oprocha2013b}]
For every $n\geq 2$, there exists a minimal subshift on $n$ symbols such that
it contains a perfect weakly mixing set of order $n$ but no non-trivial weakly mixing set of order $n+1$.
\end{thm}

Recall that it was shown in~\cite{Akin2003} that if $(X,T)$ is weakly mixing, then for every $x\in X$, the
proximal cell of $x$ is residual in $X$. In~\cite{Oprocha2013Z} Oprocha and Zhang proved that for every
closed weakly mixing set $A$ and every $x \in A$, the proximal cell of $x$ in $A$ is residual in $A$.
In \cite{Li2013c} Li, Oprocha and Zhang showed that the same is true if we consider proximal tuples instead of pairs.
First recall that an $n$-tuple $(x_1,x_2,\dotsc,x_n)\in X^n$  is \emph{proximal} if
\[\liminf_{k\to\infty}\max_{1\leq i<j\leq n} d(T^k(x_i), T^k(x_j)) = 0.\]
For $x\in  X$, define the $n$-proximal cell of $x$ as
\[P_n(x)=\{(x_1,\dotsc,x_{n-1})\in X^{n-1}\colon  (x, x_1,\dotsc,x_{n-1}) \text{ is proximal}\}.\]

\begin{thm}[\cite{Li2013c}]
Let $(X,T)$ be a dynamical system and $A\subset X$ be a weakly mixing set.
Then for every $x \in A$ and $n \geq  2$, the set $P_n(x)\cap A^{n-1}$  is residual in $A^{n-1}$.
\end{thm}

In fact, they proved even more as presented in the following theorem, where
$LY^\delta_n(x)$ is the $n$-scrambled cell of $x$ with modular $\delta>0$.

\begin{thm}[\cite{Li2013c}]
Let $(X,T)$ be a dynamical system and $A\subset X$ be a weakly mixing set.
Then for every $n\geq 2$, there exists some $\delta>0$ such that for every $x\in A$,
the set $LY_n^\delta(x)\cap A^{n-1}$  is residual in $A^{n-1}$.
\end{thm}

The following result shows that, when we look only at separation of trajectories of tuples,
weak mixing of order $2$ is enough to obtain rich structure of such points.

\begin{thm}[\cite{Li2013c}]
Let $(X,T)$ be a dynamical system and $A\subset X$ a weakly mixing set of order $2$. Then
$A$ is a sensitive set in $\bigl(\overline{\mathrm{Orb}(A,T)},T\bigr)$,
where $\mathrm{Orb}(A,T)=\{T^nx\colon n\geq 0,x\in A\}$.
In particular, the system $\bigl(\overline{\mathrm{Orb}(A,T)},T\bigr)$ is $n$-sensitive for every $n\geq 2$.
\end{thm}

Recall that we have given the definition of independent sets in Definition~\ref{def:ind-sets}.
Huang, Li and Ye in~\cite{Huang2012} studied independent sets via Furstenberg families.
In particular, they showed the following connection between
weak mixing and independent sets of open sets.
\begin{thm} \label{thm:wm-ind-ip}
Let $(X,T)$ be a dynamical system. Then the following conditions are equivalent:
\begin{enumerate}
  \item $(X,T)$ is weakly mixing;
  \item for any two non-empty open subsets $U_1,U_2$ of $X$,
  $(U_1,U_2)$ has an infinite independent set;
  \item for any $n\in\N$ and non-empty open subsets $U_1,U_2,\dotsc,U_n$ of $X$,
  $(U_1,U_2,\dotsc,U_n)$ has an IP-independent set.
\end{enumerate}
\end{thm}

In the spirit of~\cite{Huang2012} we introduce a local version of independence sets as follows.

\begin{defn}
Let $(X,T)$ be a dynamical system and $A\subset X$.
Let $U_1,U_2,\dotsc,U_n$ be open subsets of $X$ intersecting $A$.
We say that a non-empty subset $F$ of $\Z$ is an
\emph{independence set for $(U_1,U_2,\dotsc,U_n)$ with respect to $A$},
if for every non-empty finite subset $J\subset F$, and $s\in\{1,2,\dotsc,n\}^J$,
\[\bigcap_{j\in J} T^{-j}(U_{s(j)})\]
is a non-empty open subset of $X$ intersecting $A$.
\end{defn}

Now we can employ introduced notion, to state a theorem analogous to Theorem~\ref{thm:wm-ind-ip}.

\begin{thm}[\cite{Li2013c}] \label{thm:WM-ip}
Let  $(X,T)$ be a dynamical system and $A\subset X$ a closed set.
Then the following conditions are equivalent:
\begin{enumerate}
  \item\label{enum:Awm} $A$ is a weakly mixing set;
  \item\label{enum:Ind2} for every $n\geq 2$ and every open subsets $U_1,U_2,\dotsc,U_n$ of $X$ intersecting $A$,
  there exists $t\in\N$ such that $\{0,t\}$ is an independence set for $(U_1,U_2,\dotsc,U_n)$ with respect to $A$;
  \item\label{enum:IndIP} for every $n\geq 2$ and every open subsets $U_1,U_2,\dotsc,U_n$ of $X$ intersecting $A$,
  there exists a sequence $\{t_j\}_{j=1}^\infty$ in $\N$   such that
  $\{0\}\cup FS(\{t_j\}_{j=1}^\infty)$ is an independence set for $(U_1,U_2,\dotsc,U_n)$ with respect to $A$.
\end{enumerate}
\end{thm}

\subsection{Weakly mixing sets via Furstenberg families}

In 2004, Huang, Shao and Ye generalized Theorem~\ref{thm:XY1991} to $\F$-mixing systems.
Let $\F$ be a Furstenberg family. Recall that a dynamical system $(X,T)$ is called \emph{$\F$-mixing}
if it is weakly mixing and for any two non-empty open subsets $U$ and $V$ of $X$,
$N(U,V)\in \F$.

\begin{thm}[\cite{Huang2004a}] \label{thm:f-mixing}
Let $(X,T)$ be a non-trivial dynamical system and $\F$ be a Furstenberg family.
Then $(X,T)$ is $\F$-mixing if and only if for any $S\in\kappa\F$,
there exists a dense Mycielski subset $C$ of $X$
satisfying for any subset $D$ of $C$ and any continuous map $f\colon D\to X$, there exists an
increasing sequence $\{q_i\}$ in $S$ such that
$\lim\limits_{i\to\infty} T^{q_i}x=f(x)$ for all $x\in D$.
\end{thm}

It is natural that weakly mixing sets can be also generalized via Furstenberg families.
\begin{defn}
Let $(X,T)$ be a dynamical system and $\F$ be a Furstenberg family.
Suppose that $A$ is a closed subset of
$X$ with at least two points.
The set $A$ is said to be \emph{$\F$-mixing} if for any $k\in\N$,
any open subsets $U_1,U_2,\dots,U_k$, $V_1,V_2,\dotsc,V_k$ of $X$ intersecting $A$,
\[\bigcap_{i=1}^k N(U_i\cap A,V_i)\in \F.\]
\end{defn}

Inspired by the proof of Theorem \ref{thm:f-mixing},
we have the following characterization of $\F$-mixing sets.

\begin{thm}[\cite{Li2014}]\label{thm:F-Ming-set}
Let $(X,T)$ be a dynamical system and $\F$ be a Furstenberg family.
Suppose that $A$ is a closed subset of
$X$ with at least two points.
Then $A$ is an $\F$-mixing set if and only if for every $S\in\kappa\F$ (the dual family of $\F$)
there are Cantor subsets $C_1\subset C_2\subset \dotsc$ of $A$ such that
\begin{enumerate}
\item[(i)] $K=\bigcup_{n=1}^\infty C_n$ is dense in $A$;
\item[(ii)] for any $n\in\N$ and any continuous function $g\colon C_n\to A$
there exists a subsequence $\{q_i\}$ of $S$ such that
$\lim_{i\to\infty} T^{q_i}(x)=g(x)$ uniformly on $x\in C_n$;
\item[(iii)] for any subset $E$ of $K$ and
any continuous map $g \colon E\to A$ there exists a subsequence $\{q_i\}$ of $S$ such that
$\lim_{i\to\infty} T^{q_i}(x)=g(x)$ for every $x\in E$.
\end{enumerate}
\end{thm}

It is shown in~\cite{Li2014} that two classes of important dynamical systems have
weakly mixing sets via proper Furstenberg families.

Let $F$ be a subset of $\Z$. The \emph{upper Banach density of $F$} is defined by
\[BD^*(F)=\limsup_{|I|\to\infty}\frac{|F\cap I|}{|I|},\]
where $I$ is taken over all non-empty finite intervals of $\Z$.
The family of sets with positive upper Banach density is denoted by
$\F_{pubd}=\{F\subset \Z\colon BD^*(F)>0 \}$.
We say that $F$ is \emph{piecewise syndetic} if it is the intersection of a thick set and a syndetic set.
The family of piecewise syndetic sets is denoted by $\F_{ps}$.
\begin{thm}[\cite{Li2014}]
Let $(X,T)$ be a dynamical system.
\begin{enumerate}
  \item If $(X,T)$ has positive topological entropy, then it has some $\F_{pubd}$-mixing sets.
  \item If $(X,T)$ is a non-PI minimal system, then it has some $\F_{ps}$-mixing sets.
\end{enumerate}
\end{thm}

\section{Chaos in the induced spaces}

A dynamical system $(X,T)$ induces two natural systems,
one is $(K(X),T_K)$ on the hyperspace $K(X)$ consisting of all closed non-empty subsets of
$X$ with the Hausdorff metric, and the other one is $(M(X),T_M)$ on
the probability space $M(X)$ consisting of all Borel probability
measures with the weak*-topology.
Bauer and Sigmund \cite{Bauer1975}
first gave a systematic investigation on the connection of dynamical
properties among $(X,T)$, $(K(X),T_K)$ and $(M(X),T_M)$. It was
proved that $(X,T)$ is weakly mixing (resp. mildly mixing, strongly
mixing) if and only if $(K(X),T_K)$ (resp. $(M(X),T_M)$) has the
same property.
This leads to a natural question:

\begin{prob}\label{pp6}
If one of the dynamical systems $(X,T)$, $(K(X),T_K)$ and $(M(X),T_M)$ is chaotic in some sense,
how about the other two systems?
\end{prob}

This question attracted a lot of attention, see, e.g.,
\cite{Roman-Flores2003,Banks2005,Guirao2009} and references therein, and many partial answers were obtained.
We first show that when the induced system is weakly mixing.
\begin{thm}[\cite{Bauer1975,Banks2005}]
Let $(X,T)$ be a dynamical system.
Then $(X,T)$ is weakly mixing if and only if $(K(X),T_K)$ is weakly mixing if and only if $(K(X),T_K)$ is transitive.
\end{thm}

\begin{thm}[\cite{Bauer1975,LiJ2013}]
Let $(X,T)$ be a dynamical system.
Then $(X,T)$ is weakly mixing if and only if $(M(X),T_M)$ is weakly mixing if and only if $(M(X),T_M)$ is transitive.
\end{thm}

In \cite{Glasner1995},  Glasner and Weiss studied the topological entropy of $(K(X),T_K)$ and $(M(X),T_M)$.
They proved that
\begin{thm}
Let $(X,T)$ be a dynamical system.
\begin{enumerate}
  \item The topological entropy of $(X,T)$ is zero if and only if the one of $(M(X),T_M)$ is also zero,
  and the topological entropy of $(X,T)$ is positive if and only if the one of $(M(X),T_M)$ is infinite.
  \item If the topological entropy of $(X,T)$ is positive, then the topological entropy of $(K(X),T_K)$  is infinite,
  while there exists a minimal system $(X,T)$ of zero topological entropy and
  $(K(X),T_K)$ with positive topological entropy.
\end{enumerate}
\end{thm}

To show that when the dynamical system on the hyperspace is Devaney chaotic,
we need to introduce some concepts firstly.
We say that a dynamical system $(X, T)$ has \emph{dense small
periodic sets} \cite{Huang2005} if for any non-empty open subset $U$ of
$X$ there exists a non-empty closed subset $Y$ of $U$ and $k \in \mathbb{N}$
such that $T^kY \subset Y$. Clearly, if a dynamical system has a dense set of periodic points,
then it also has dense small periodic sets.
The dynamical system $(X, T)$ is called an \emph{HY-system} if it is totally
transitive and has dense small periodic sets.
Note that there exists an HY-system without periodic points (see~\cite[Example~3.7]{Huang2005}).
Recently, Li showed in~\cite{Li2014d} that $(K(X),T_K)$ is Devaney chaotic is equivalent to
the origin system $(X,T)$ is an HY-system.

\begin{thm}[\cite{Li2014d}]
Let $(X,T)$ be a dynamical system with $X$ being
infinite. Then the following conditions are equivalent:
\begin{enumerate}
  \item $(K(X), T_K)$ is Devaney chaotic;
  \item $(K(X),T_K)$ is an HY-system;
  \item $(X,T)$ is an HY-system.
\end{enumerate}
\end{thm}

In order to characterize Devaney chaos on the space of probability measures, we need
a notion of an almost HY-system.
We say that $(X,T)$ has \emph{almost dense
periodic sets} if for each non-empty  open subset $U \subset X$ and
$\epsilon>0$, there are $k \in \N$ and $\mu \in M(X)$ with
$T_M^k\mu=\mu$ such that $\mu(U^c)<\epsilon$, where $U^c=\{x \in X:
x \notin U\}$. We say that $(X, T)$ is an \emph{almost HY-system} if
it is totally transitive and has almost dense periodic sets.

\begin{thm}[\cite{LiJ2013}]
Let $(X,T)$ be a dynamical system with $X$ being
infinite. Then the following conditions are equivalent:
\begin{enumerate}
  \item $(M(X), T_M)$ is Devaney chaotic;
  \item $(M(X),T_M)$ is an almost HY-system;
  \item $(X,T)$ is an almost HY-system.
\end{enumerate}
\end{thm}
It is clear that every HY-system is also an almost HY-system.
There is a non-trivial minimal weakly mixing almost-HY-system (see \cite[Theorem 4.11]{LiJ2013}),
which is not an HY-system, since every minimal HY-system is trivial.

Recall that a dynamical system $(X,T)$ is \emph{proximal} if any pair $(x,y)\in X^2$ is proximal.
The following proposition shows that
if $(K(X),T_K)$ is proximal then $(X,T)$ is ``almost'' trivial.

\begin{prop}[\cite{Li2014a}]
Let $(X,T)$ be a dynamical system. Then the following conditions are equivalent:
\begin{enumerate}
  \item $(K(X),T_K)$ is proximal;
  \item $\bigcap_{n=1}^\infty T^n X$ is a singleton;
  \item $X$ is a uniformly proximal set.
\end{enumerate}
\end{prop}

If $(M(X),T_M)$ is proximal, then $(X,T)$ is called \emph{strongly proximal}~\cite{Glasner1976}.
Note that if $(X,T)$ is strongly proximal, then it is proximal,
since $(X,T)$ can be regarded as a subsystem of the proximal system $(M(X),T_M)$.
We have the following characterization of strongly proximal systems.
\begin{thm}[\cite{Li2014a}]
Let $(X,T)$ be a dynamical system. Then the following conditions are equivalent:
\begin{enumerate}
  \item $(X,T)$ is strongly proximal;
  \item $(X,T)$ is proximal and unique ergodic;
  \item every pair $(x,y)\in X^2$ is Banach proximal.
\end{enumerate}
\end{thm}

\bigskip
\noindent\textbf{Acknowledgement: }
The authors would like to thank Wen Huang, Jie Li and Guohua Zhang for the careful reading,
and thank the anonymous reviewer for his/her valuable comments and suggestions.
The first author is partially supported by NNSF of China (11401362, 11471125)
and NSF of Guangdong province (S2013040014084).
The second author is partially supported by NNSF of China (11371339,11431012).

\bibliographystyle{amsplain}
\providecommand{\bysame}{\leavevmode\hbox to3em{\hrulefill}\thinspace}
\providecommand{\MR}{\relax\ifhmode\unskip\space\fi MR }
% \MRhref is called by the amsart/book/proc definition of \MR.
\providecommand{\MRhref}[2]{%
  \href{http://www.ams.org/mathscinet-getitem?mr=#1}{#2}
}
\providecommand{\href}[2]{#2}

\end{document}